\documentclass[hidelinks,onefignum,onetabnum]{siamart220329}

\usepackage{lipsum}
\usepackage{amsfonts}
\usepackage{graphicx}
\usepackage{epstopdf}
\usepackage{algorithmic}
\ifpdf
  \DeclareGraphicsExtensions{.eps,.pdf,.png,.jpg}
\else
  \DeclareGraphicsExtensions{.eps}
\fi

\newsiamremark{remark}{Remark}

\newsiamremark{hypothesis}{Hypothesis}
\crefname{hypothesis}{Hypothesis}{Hypotheses}
\newsiamthm{claim}{Claim}

\headers{FHT under wavelet basis}{X. Tang, L. Ying}

\title{Wavelet-based density sketching with functional hierarchical tensor}

\author{
  Xun Tang\thanks{Corresponding author. Institute for Computational and Mathematical Engineering (ICME), Stanford University, Stanford, CA 94305, USA. 
  (\email{xuntang@stanford.edu})}
  \and
  Lexing Ying\thanks{Department of Mathematics and Institute for Computational and Mathematical Engineering (ICME), Stanford University, Stanford, CA 94305, USA. 
  (\email{lexing@stanford.edu})}
  \funding{X.T. and L.Y. are supported by AFOSR MURI award FA9550-24-1-0254.}
  }

\usepackage{amsopn}

\makeatletter
\newcommand*{\addFileDependency}[1]{%
  \typeout{(#1)}%
  \@addtofilelist{#1}%
  \IfFileExists{#1}{}{\typeout{No file #1.}}%
}
\makeatother

\usepackage{amsmath,amsfonts,bm}

\def\eqref#1{equation~\ref{#1}}

\def\ceil#1{\lceil #1 \rceil}

\def\1{\bm{1}}

\DeclareMathAlphabet{\mathsfit}{\encodingdefault}{\sfdefault}{m}{sl}
\SetMathAlphabet{\mathsfit}{bold}{\encodingdefault}{\sfdefault}{bx}{n}

\newcommand{\R}{\mathbb{R}}

\usepackage{amsmath,amssymb,amsfonts}
\usepackage{latexsym}
\usepackage{indentfirst}
\usepackage{graphicx}
\usepackage{placeins}
\usepackage{enumerate}
\usepackage{booktabs}
\usepackage{algorithm}
\usepackage{algorithmic}
\usepackage{multirow}
\usepackage{optidef}
\usepackage{hyperref}
\usepackage{graphicx,color}
\usepackage{diagbox}
\usepackage{subcaption}
\usepackage[normalem]{ulem}

\usepackage{cleveref}

\ifpdf
\hypersetup{
  pdftitle={Multi-scale-FTN},
  pdfauthor={X. Tang, and L. Ying}
}
\fi

\begin{document}

\maketitle

\begin{abstract}
    We introduce the functional hierarchical tensor under a wavelet basis (FHT-W) ansatz for high-dimensional density estimation in lattice models. Recently, the functional tensor network has emerged as a suitable candidate for density estimation due to its ability to calculate the normalization constant exactly, a defining feature not enjoyed by neural network alternatives such as energy-based models or diffusion models.
    While current functional tensor network models show good performance for lattice models with weak or moderate couplings, we show that they face significant model capacity constraints when applied to lattice models with strong coupling. To address this issue, this work proposes to perform density estimation on the lattice model under a wavelet transformation.
    Motivated by the literature on scale separation, we perform iterative wavelet coarsening to separate the lattice model into different scales. Based on this multiscale structure, we design a new functional hierarchical tensor ansatz using a hierarchical tree topology, whereby information on the finer scale is further away from the root node of the tree. 
    Our experiments show that the numerical rank of typical lattice models is significantly lower under appropriate wavelet transformation.
    Furthermore, we show that our proposed model allows one to model challenging Gaussian field models and Ginzburg-Landau models.
\end{abstract}

\begin{keyword}
    High-dimensional density estimation; Functional tensor network; Wavelet transform.
\end{keyword}

\begin{MSCcodes}
    15A69, 65F99, 42C40
\end{MSCcodes}

\section{Introduction}\label{sec: introduction}

Estimating high-dimensional density is one of the most important tasks in applied sciences. It is used in various fields of data science and scientific computing.
During density estimation, one is given a collection of samples from a target distribution \(p \colon \R^{d} \to \R\), and the goal is to use the samples to obtain an accurate approximation of \(p\) from a chosen function class.

The functional tensor network (FTN) has recently emerged as an appealing function class to model probability distributions. The key feature of an FTN model is that its normalization constant calculation can be done via an efficient tensor contraction diagram, which is especially important for obtaining normalized density functions. For FTN models based on line graph \cite{oseledets2010tt,bigoni2016spectral} and the binary tree graph \cite{hackbusch2009new,tang2024solving}, calculating the normalization constant has an \(O(d)\) scaling, where \(d\) stands for the dimension of the system. In contrast, performing an accurate numerical integration of parametric functions typically suffers from an exponential scaling in \(d\).

This work focuses on distributions coming from lattice models in 1D and 2D, for which FTN models are known to be suitable numerical candidates for density estimation tasks \cite{chen2023combining,tang2024solving}. Lattice models in 1D and 2D often come from discretizing a stochastic process. For example, let \(\{X_t \colon t \in [0, 1]\}\) be a Gaussian process, e.g., the Ornstein–Uhlenbeck process \cite{williams2006gaussian}, and then \(x = (X_{i/d})_{i = 1}^{d}\) constitutes a 1D lattice model. Similar constructions can be done in 2D or higher physical dimensions. Overall, discretization of smaller domains leads to lattice models with stronger site-wise coupling.

This work is primarily motivated by difficulties in using existing tensor network models to tackle lattice representations with strong coupling. 
In practice, a 1D lattice model is often of the form \[p(x_1, \ldots, x_d) = \exp\left(-\frac{\alpha}{2}\sum_{|i - j| = 1}(x_i - x_j)^2 - \sum_{i = 1}^{d}V_i(x_i)\right),\] 
and the case where \(\alpha\) is large leads to situations where one needs to increase the rank of the FTN models to represent \(p\) faithfully. To give \(p\) a more efficient low-rank approximation, we propose to use the wavelet transformation, which has found great uses in image analysis and signal processing \cite{daubechies1992ten,mallat1999wavelet}. A wavelet transformation can localize information both in frequency and in space, which allows one to model \(p\) with a smaller rank requirement. The wavelet transformation is orthogonal, and therefore, it is an invertible and measure-preserving transformation.

For a 1D lattice model with \(d = 2^{L}\), successive wavelet coarse-graining transforms the physical \(x\) space into the wavelet coordinate \((c_{k, l})_{k, l}\). The \(l = -1, \ldots, L-1\) index is the resolution scale index, with larger \(l\) standing for information on a finer scale. The \(k = 1, \ldots, 2^{\max(l, 0)}\) index is the location index at scale \(l\). The recent work in \cite{marchand2023multiscale} shows that a multiscale representation with successive wavelet-based coarse-graining can well approximate 1D and 2D lattice models. Moreover, for practical lattice models, the analysis in \cite{marchand2023multiscale} shows that the interaction between the variables \(c_{k, l}\) and \(c_{k', l'}\) is only significant when \((k, l)\) is close to \((k', l')\) in both indices. The localized interaction is an instance of the well-celebrated scale separation phenomenon in renormalization group (RG) theory \cite{wilson1971renormalization,wilson1972critical,kadanoff1976variational}. 

The proposed iterated wavelet transform, commonly referred to as wavelet multiresolution approximation \cite{mallat1999wavelet}, calls for a novel design of FTN topology. The transformed variable \(c = (c_{k, l})_{k, l}\) has a natural hierarchical structure suitable for a binary tree graph reminiscent of the functional hierarchical tensor model \cite{tang2024solving}. However, the current FHT model places the external bond on the leaves nodes of a binary tree, which is less suitable. Therefore, we introduce a new functional tensor network structure where the variables correspond to nodes on a modified binary tree graph. We refer to the proposed function class as the functional hierarchical tensor under a wavelet basis (FHT-W). We introduce the density algorithm for FHT-W, and for future reference, we go through the subroutine of density estimation for FTN models under arbitrary tree structures. Our numerical experiments show that the FHT-W ansatz can successfully model high-dimensional 1D and 2D lattice models.

\subsection{Motivating example}\label{sec: motivating example}
For FTN models, the most important model capacity parameter is the maximal internal bond dimension \(r\). For a \(d\)-dimensional FHT model, the parameter size is \(O(dr^3)\), and for functional tensor train (FTT) \cite{bigoni2016spectral} the parameter size is \(O(dr^2)\). A coordinate system with a better low-rank structure leads to a significant practical speedup both in storage complexity and computational complexity. In this light, an FTN ansatz's ability to capture a 1D or 2D lattice model largely depends on whether it can represent the model with a reasonably small rank \(r\). 

We illustrate how coupling strength influences rank through a simple 2D distribution function \(p(x_1, x_2) = e^{-\frac{x_1^2}{2} - \frac{x_2^2}{2} - \frac{\alpha}{2}(x_1 - x_2)^2}\) with \(\alpha > 0\). In this setting, using an FTT or FHT approximation for \(p\) is equivalent to the task of finding functions \(h_i(x_1), g_i(x_2)\) for \(i =1,\ldots, r\) so that \(p(x_1, x_2) \approx \sum_{i=1}^{r}h_i(x_1)g_i(x_2)\). 
When \(\alpha \) is small, the system has weak coupling between \(x_1\) and \(x_2\), and it is standard to show that \(p\) admits a low-rank approximation, e.g. by a Chebyshev approximation of the exponential function \cite{demanet2010chebyshev}. However, as \(\alpha\) increases, it becomes increasingly difficult to approximate \(p\) in a low-rank representation unless one increases \(r\) with \(\alpha\).
Write \(X = (X_1, X_2) \sim p\). One sees that \(X_1 \sim N(0, \frac{1+\alpha}{1+2\alpha})\) and \(X_2 | X_1 = c \sim N(c, \frac{1}{1+2\alpha})\). Therefore, as \(\alpha \to \infty\), one has \(X\) converging in distribution to \(W = (W_1, W_2)\), where \(W_1 \sim N(0, 1/2)\) and \(W_2 = W_1\). One can see that the distribution of \(W\) is singular and does not admit a low-rank approximation.

While the large rank issue occurring for \(p\) in the limiting case might seem discouraging for tensor network models, one can address the issue in this case by a coordinate transformation.
Writing \((\xi_1, \xi_2) = (\frac{x_1 - x_2}{\sqrt{2}}, \frac{x_1 + x_2}{\sqrt{2}})\), one sees that \(p(\xi_1, \xi_2) = \exp(-(\alpha + \frac{1}{2}) \xi_1^2)\exp( - \frac{\xi_2^2}{2})\) is exactly of rank \(r = 1\). Therefore, \(p\) in the \((\xi_1, \xi_2)\) coordinate always admits a low-rank representation regardless of what values of \(\alpha\) one chooses.

An extension of the two-dimensional example is the Ornstein–Uhlenbeck (O-U) model, for which one has \(p(x_1, \ldots, x_d) = \exp\left(-\frac{\alpha}{2}\sum_{i \sim j}(x_i - x_j)^2 - \frac{\beta}{2}\sum_{i = 1}^{d}x_i^2\right)\), where \(i \sim j\) if \(i - j = 1 \mod {d}\). One sees that the precision matrix of \(p\) is a circulant matrix, and so the model \(p(\xi)\) is separable if one takes \(\xi\) to be the discrete Fourier transform (DFT) of \(x\).

By the illustrated examples, one can see that performing an FTN density estimation of O-U models is practically much more efficient after a Fourier transformation. The Fourier transform is generally not sufficient for a low-rank parameterization, as DFT transforms local nonlinear terms into global interaction terms in the frequency space. In contrast, the wavelet transformation is known to give lattice models a better low-rank structure \cite{marchand2023multiscale}. In \Cref{sec: numerical rank}, we exhibit strong numerical evidence that the probability distribution under a wavelet-transformed coordinate leads to a significantly smaller numerical rank to capture the target distribution. In fact, our experiments in \Cref{sec: numerical rank} show that the numerical rank for practical 2D lattice models can be reduced by up to a factor of five if one uses the wavelet basis.

\subsection{Discussion on high-dimensional density estimation}

In addition to the key role of density estimation in generative modeling \cite{goodfellow2014generative}, the density estimation subroutine can often aid one in solving high-dimensional partial differential equations (PDE). 
For the Fokker-Planck equation, which models the evolution of probability densities of a particle undergoing advection and diffusion, one can obtain an approximate PDE solution by using density estimation in a particle-based framework \cite{tang2024solving}. In the Kolmogorov backward equation, which models the evolution of the conditional expectation function, one can use density estimation subroutines to construct the solution operator, which gives the PDE solution for arbitrary terminal conditions \cite{tang2024solvinga}.

Density estimation is a difficult task. Given samples from a high-dimensional distribution \(p\), it is well-known that nonparametric density estimation requires an exponential number of samples in the dimension \(d\) \cite{wasserman2006all}. Traditional parametric function classes such as tree-based graphical models \cite{jordan1999learning} do admit efficient density estimation through maximum likelihood estimation, but its representation power is limited. While neural networks provide several new function classes for density estimation, a neural network often represents an unnormalized likelihood function, as is true for energy-based models \cite{hinton2002training,lecun2006tutorial} and diffusion models \cite{song2021maximum}. Therefore, a key element in density estimation research lies in designing new parametric function classes that possess the ability to produce a normalized likelihood function.

When one applies density estimation to lattice models, the goal is to model a probability distribution under moderate computational resources. Gaussian field models can be well-approximated by fitting a multivariate normal distribution with maximum likelihood. However, it is often true that the lattice model is multi-modal and is hence ill-suited to be modeled by normal distributions. An example of a lattice model with multimodality is the Ginzburg-Landau (G-L) model used for studying the phenomenological theory of superconductivity \cite{ginzburg2009theory, hoffmann2012ginzburg, hohenberg2015introduction,weinan2004minimum}. In particular, the Ginzburg-Landau model contains a double-well term, which makes a multivariate normal distribution unsuitable to fit such models.

In the 2D G-L model, one can intuitively think of a particle taking the form of a field \(x(a) \colon [0,1]^2 \to \R\). To obtain a lattice model, one discretizes the unit square \([0,1]^2\) into a grid of \(d = m^2\) points $\{(ih,jh)\}$ for \(h = \frac{1}{m+1}\) and $1\le i,j\le m$. The discretization of the field at this grid is the $d$-dimensional vector $x=(x_{(i,j)})_{1\le i,j \le m}$, where \(x_{(i, j)} = x(ih, jh)\). The field dynamics under this discretization is the overdamped Langevin dynamics under the following potential function
\begin{equation}\label{eqn: 2D GZ model}
V(x) = V(x_{(1,1)}, \ldots ,x_{(m,m)}) := \frac{\lambda}{2} \sum_{v \sim w}\left(\frac{x_{v} - x_{w}}{h}\right)^2 +  \frac{1}{4\lambda} \sum_{v}\left(1 - x_v^2\right)^2,
\end{equation}
where $v$ and $w$ are Cartesian grid points and \(v \sim w\) if and only if they are adjacent. One can see that the term \((1 - x_v)^2\) has two favorable configurations of \(x_v = 1\) and \(x_v = -1\), which is the double-well term that makes multivariate normal distributions ill-suited to approximate G-L models. Moreover, a finer spatial resolution increases the \(1/h\) term in \Cref{eqn: 2D GZ model}, which leads to stronger site-wise coupling.

The work in \cite{tang2024solving} shows that FHT works well for 1D-3D Ginzburg-Landau models with moderate coupling. However, when the G-L model has a strong coupling term and a weak double-well term, one can see that the necessary rank to approximate the G-L model can be quite large. Therefore, the proposed FHT-W ansatz is a reasonable way to utilize the efficiency of the FTN models while successfully dealing with strong coupling situations in the G-L model.

For general distributions, it is less clear how one can apply a coordinate transformation to provide a good low-rank structure. However, given the inherent difficulty of density estimation tasks, FTN models remain one of the most expressive existing function classes to generate un-normalized densities, and combining FTN with coordinate transformation strategies can strongly improve the expressive power of existing FTN models \cite{dektor2023tensor,ren2023high}.

\subsection{Related work}\label{sec: related work}

\paragraph{Neural network models in density estimation}
Neural networks have been used in density estimation-related tasks within generative modeling and variational inference. Examples include restricted Boltzmann machines (RBMs)~\cite{hinton2010practical, salakhutdinov2010efficient}, energy-based models (EBMs)~\cite{hinton2002training, lecun2006tutorial, gutmann2010noise}, variational auto-encoders (VAEs)~\cite{doersch2016tutorial, kingma2019introduction}, generative adversarial networks (GANs)~\cite{goodfellow2014generative}, normalizing flows~\cite{tabak2010density, tabak2013family, rezende2015variational, papamakarios2021normalizing}, diffusion and flow-based models~\cite{sohl2015deep,zhang2018monge,song2019generative,song2020score, song2021maximum, ho2020denoising, albergo2022building, liu2022flow, lipman2022flow, albergo2023stochastic}. Among these methods, only normalizing flow can produce a normalized probability density. We remark that normalizing flow uses a neural network to parameterize a reversible flow map, which typically leads to reduced expressivity compared to other approaches. Moreover, the back-propagation of the Jacobian term during the maximum likelihood training is quite difficult compared to diffusion models.

\paragraph{Functional tensor network methods for density estimation}
The functional tensor network has been extensively explored as an ansatz for density estimation. Density estimation with functional tensor networks has been done in functional tensor trains \cite{novikov2021tensor,hur2023generative,chen2023combining} and functional hierarchical tensors \cite{tang2024solving,tang2024solvinga,sheng2025numerical}. 
The aforementioned development of density estimation subroutines in the continuous case is heavily inspired by the discrete counterpart in the case of tensor train models and hierarchical tensor models \cite{han2018unsupervised,cheng2019tree, glasser2019expressive, bradley2020modeling, grelier2022learning, tang2023generative, ren2023high, khoo2023tensorizing,peng2023generative}.

\paragraph{Wavelet transformation in lattice models}
Wavelet analysis has proven to be a critical tool in understanding lattice models through its close connection with renormalization groups \cite{battle1999wavelets,altaisky2016unifying}. In \cite{evenbly2016entanglement}, a multiscale entanglement
renormalization ansatz is constructed for the 1D lattice system of fermions through the Daubechies D4 wavelets, which results in a remarkably accurate approximation of the ground state of the critical Ising model, and a subsequent work in \cite{haegeman2018rigorous} considers the case of wavelet-based ground state approximation in 2D lattice fermionic systems. For continuous lattice models considered in this work, one of the main theoretical foundations for using wavelets for the Ginzburg-Landau model is that the wavelet transformation can nearly diagonalize the Laplacian term \cite{beylkin1992representation,meyer1992wavelets}. In \cite{marchand2023multiscale}, it is further shown that sufficiently high-order wavelet methods can localize the double well term in the G-L model.

\paragraph{Wavelet transformation in generative modeling}
Wavelet has been instrumental in generative modeling, especially in the setting of image synthesis. Wavelet has been used in normalizing flow \cite{yu2020wavelet,dai2024multiscale}, diffusion models \cite{guth2023conditionally,kadkhodaie2023learning}, energy-based models \cite{marchand2023multiscale}, and auto-encoders \cite{chen2018learning}. Other related works include \cite{ zhang2021multiscale,gal2021swagan,ho2022cascaded,saharia2022image}.
Among these, this work bears the most resemblance to the work in \cite{marchand2023multiscale}, which performs generative modeling under a wavelet representation using a multiscale energy-based model. Compared to \cite{marchand2023multiscale}, while this work likewise performs a wavelet transformation on the input data, this work uses a tensor network model so that one can perform density estimation under a wavelet representation.

\subsection{Outline}
This work is organized as follows. \Cref{sec: density estimation} details the density estimation subroutine for FTN models under general tree topology. \Cref{sec: main formulation} covers the proposed functional hierarchical tensor used for wavelet-transformed samples. \Cref{sec: numerical rank} presents numerical experiments that show wavelet transformation decreases the numerical rank of lattice models. \Cref{sec: numerical experiment} shows the performance of the proposed ansatz on practical 1D and 2D lattice models with strong coupling.

\section{Density estimation under general tree topology}\label{sec: density estimation}

This section details the sketching-based density estimation subroutine for functional tensor network models under general tree topology. The FHT-W ansatz to be introduced in \Cref{sec: main formulation} is a special case that takes on a specific tree topology. This generalized setting allows for a cleaner presentation and accommodates other possible tree topologies in potential future works. 
The introduced subroutine is also a generalized version of the density estimation routine considered in \cite{hur2023generative,chen2023combining, tang2024solving}.

\paragraph{Notation}
For notational compactness, we introduce several shorthand notations for simple derivation. For \(n \in \mathbb{N}\), let \([n] := \{1,\ldots, n\}\). For an index set \(S \subset [d]\), we let \(x_{S}\) stand for the subvector with entries from index set $S$. The symbol \(\Bar{S}\) stands for the set-theoretic complement \(\Bar{S} = [d] - S\). For a \(d\)-tensor \(D \colon \prod_{j = 1}^{d}[n_j] \to \R\), we use \(D(i_S; i_{\bar{S}}) \colon [\prod_{j \in S}n_j] \times [\prod_{j \in \Bar{S}}n_j] \to \R \) to denote the \emph{unfolding matrix} by merging the indices in \(i_S\) to the row index and merging the indices in \(i_{\bar{S}}\) to the column index.

\paragraph{Functional tensor network}
We first introduce the general functional tensor network (FTN) ansatz.
Let \(\{\psi_{i,j}\}_{i = 1}^{n_j}\) denote a collection of orthonormal function basis over a single variable for the \(j\)-th coordinate, and let \(D \colon \prod_{j = 1}^{d} [n_j] \to \R\) be a \(d\)-tensor.
The functional tensor network is the \(d\)-dimensional function defined by the following equation:
\begin{equation}
\label{eqn: FTN forward map}
    p(x_{1}, \ldots, x_{d}) = \sum_{i_{k} \in [n_k], k \in [d]} D_{i_1,\ldots, i_{d}} \psi_{i_1, 1}(x_1)\cdots \psi_{i_{d}, d}(x_{d}) = \left<D, \, \bigotimes_{j=1}^{d} \Vec{\Psi}_{j}(x_j) \right>,
\end{equation}
where \(\Vec{\Psi}_{j}(x_j) = \left[\psi_{1,j}(x_j),\ldots, \psi_{n,j}(x_j) \right]\) is an \(n\)-vector encoding the evaluation for the \(x_j\) variable over the entire \(j\)-th functional basis. The definition is general, and one typically asserts special structures in \(D\) to ensure that calculations in \(p\) are efficient.

\paragraph{Tree structure notation}
We provide notations for a tree graph. A tree graph \(T = (V,E)\) is a connected undirected graph without cycles. An undirected edge \(\{v,v'\} \in E\) is written interchangeably as $(v, v')$ or $(v', v)$. For any \(v \in V\), define \(\mathcal{N}(v)\) to be the neighbors of \(v\). Moreover, define \(\mathcal{E}(v)\) as the set of edges incident to \(v\). For an edge \(e = (k, w)\), removing the edge \(e\) in \(E\) results in two connected components \(I_{1}, I_{2}\) with \(I_1 \cup I_2 = V\). 
For any \((w, k) \in E\), we use \(k \to w\) to denote the unique connected component among \(I_1, I_2\) which contains \(k\) as a node.

\begin{figure}
    \centering
    \subfloat[\centering]{\includegraphics[width=0.45\linewidth]{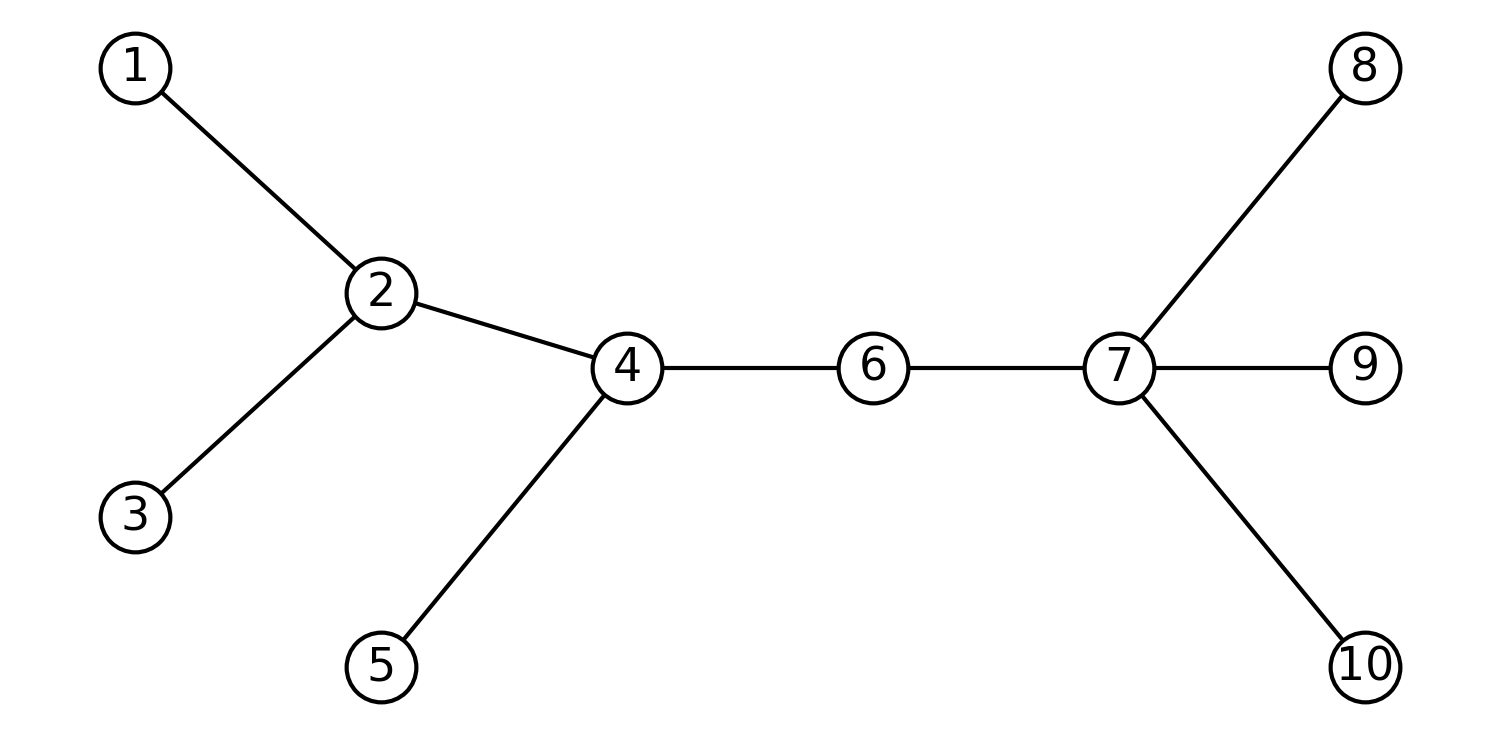}}%
    \hspace{12pt}
    \subfloat[\centering]{\includegraphics[width=0.45\linewidth]{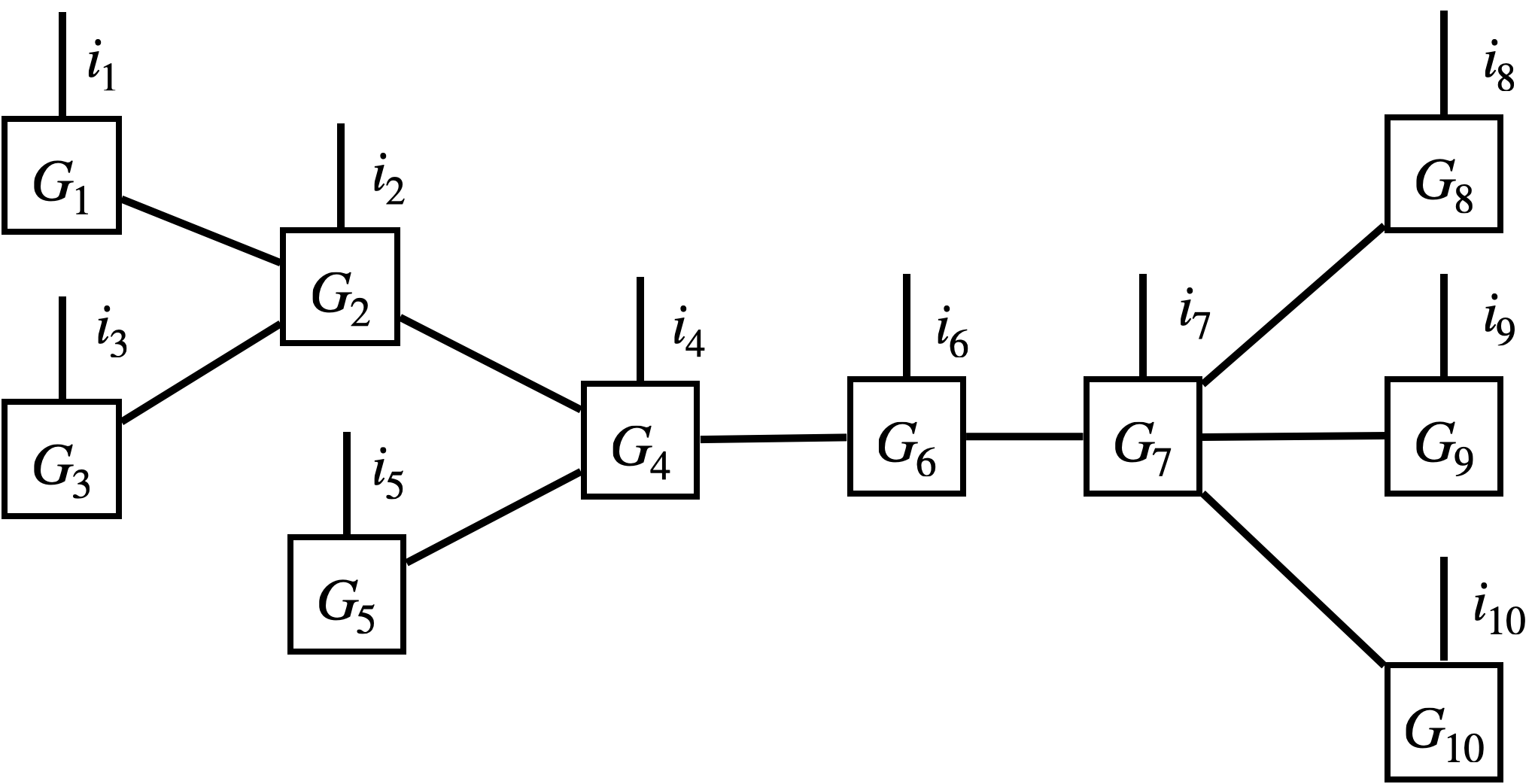}}%
    \caption{(A) A tree structure \(T = (V, E)\) with \(V = V_{\mathrm{ext}}= \{1, \ldots, 10\}\). (B) Tensor Diagram representation of a tree tensor network over \(T\).}
    \label{fig:tree+ttn_SEC_TTN}
\end{figure}

\paragraph{Tree-based FTN}
We introduce the tree-based functional tensor network ansatz to be used for density estimation. Let \(T = (V, E)\) be a tree graph and let \(V_{\mathrm{ext}}\) be the subset of vertices in \(T\) whose associated tensor component admits an external bond. To adhere to the setting of \(d\)-dimensional functions, let \(|V_{\mathrm{ext}}| = d\). The total vertex size \(\tilde{d} = |V|\) is equal to \(d\) only when the tensor network has no internal nodes. We define the tree-based functional tensor network in \Cref{def: tree-based FTN}. The tensor network is illustrated in \Cref{fig:tree+ttn_SEC_TTN}
\begin{definition}\label{def: tree-based FTN}
    (Tree-based functional tensor network)
    Suppose one has a tree structure $T = (V, E)$ and the corresponding ranks $\{r_e : e \in E\}$. The symbol \(V_{\mathrm{ext}}\) is the subset of \(V\) which contains external bonds. For simplicity, we apply a labeling of nodes so that \(V = [\tilde{d}]\) and \(V_{\mathrm{ext}} = [d]\).
    
    The \emph{tensor component} at \(k \in V\) is denoted by \(G_{k}\). Let \(\mathrm{deg}(k)\) stand for the degree of \(k\) in \(T\).
    When \(k \in V_{\mathrm{ext}}\), \(G_{k}\) is defined as a \((\mathrm{deg}(k) + 1)\)-tensor of the following shape:
    \begin{equation*}\label{eqn: core size constraint case external}
        G_{k} : [n_{k}] \times \prod_{e \in \mathcal{E}(k)}  [r_e] \rightarrow \R.
    \end{equation*}
    When \(k \not \in V_{\mathrm{ext}}\), \(G_{k}\) is defined as a \(\mathrm{deg}(k)\)-tensor of the following shape:
    \begin{equation*}\label{eqn: core size constraint case internal}
        G_{k} : \prod_{e \in \mathcal{E}(k)}  [r_e] \rightarrow \R.
    \end{equation*}
    
    A \(d\)-tensor \(D\) is said to be a \emph{tree tensor network} defined over the tree \(T\) and tensor components $\{G_k\}_{k=1}^{\tilde{d}}$ if
    \begin{equation}
        \label{eq:ttn-contraction_sec_TTN_with_internal}
        D(i_1, \ldots, i_{d}) = \sum_{\alpha_E} \prod_{k 
        =1}^{d} G_k\left(i_k, \alpha_{\mathcal{E}(k)}\right)\prod_{k = d+1}^{\tilde{d}}G_k\left( \alpha_{\mathcal{E}(k)}\right).
    \end{equation}

    Moreover, for \(j = 1, \ldots, d\), let \(\{\psi_{i,j}\}_{i = 1}^{n_j}\) denote a collection of orthonormal function basis over a single variable for the \(j\)-th coordinate. A tree-based FTN is the \(d\)-dimensional function defined by \Cref{eqn: FTN forward map} with \(D\) as the coefficient tensor. Specifically, one has
    \begin{equation}
        \label{eq: tree-based FTN forward map}
        p(x_1, \ldots, x_{d}) = \sum_{ \alpha_E} \prod_{k 
        =1}^{d} \left( \sum_{i_{k} = 1}^{n_k}G_k\left(i_k, \alpha_{\mathcal{E}(k)}\right)\psi_{i_k, k}(x_k)\right)\prod_{k = d+1}^{\tilde{d}}G_k\left( \alpha_{\mathcal{E}(k)}\right).
    \end{equation}
\end{definition}

\paragraph{Sketched linear equation}
We cover the procedure to perform density estimation over a tree-based FTN ansatz. 
The input is a collection of independent empirical samples \( \left\{ x^{(j)}:=\left(x_1^{(j)}, \ldots, x_{d}^{(j)} \right) \right\}_{j = 1}^{N}\) sampled according to an underlying distribution \(p\). In this work, the input would be samples from lattice models after the iterative wavelet transformation. We assume that \(p \colon \R^{d} \to \R\) is a tree-based FTN defined over \(D\), and \(D\) is defined over tensor components \(\{G_k\}_{k =1}^{\tilde{d}}\) as in \Cref{def: tree-based FTN}.

We shall derive the linear equation for each \(G_{k}\). Let \(e_1, \ldots, e_{\mathrm{deg}(k)}\) be all edges incident to \(k\), and let \(e_{j} = (v_j, k)\) for \(j = 1, \ldots, \mathrm{deg}(k)\). For each \(v_{j}\), we construct \(D_{v_{j} \to k}\) by contracting all tensor components \(G_{w}\) for \(w \in v_{j} \to k\). When \(k \in V_{\mathrm{ext}} = [d]\), one has the following linear equation:
\begin{equation}\label{eqn: core-determining equation}
    D(i_1, \ldots, i_d) = \sum_{a_{\mathcal{E}(k)}}G_{k}(i_k, a_{\mathcal{E}(k)}) \prod_{j = 1}^{\mathrm{deg(k)}}D_{v_j \to k}(\alpha_{e_j},i_{v_j \to k \cap [d]}),
\end{equation}
and the equation for \(k = d + 1 ,\ldots, \tilde{d}\) can be obtained from \Cref{eqn: core-determining equation} by simply omitting the \(i_{k}\) index in the right hand side of the equation. 

We first consider the case for \(k \in V_{\mathrm{ext}} = [d]\). We see that \Cref{eqn: core-determining equation} is exponential-size. Thus, we apply sketching to obtain a linear system of reasonable size.
Let \(\mathcal{E}(k) = \{e_{j} = (v_j, k)\}_{j = 1}^{\mathrm{deg(k)}}\). For \(j = 1, \ldots, \mathrm{deg}(k)\), we let \(\tilde{r}_{e_{j}} > r_{e_{j}}\) be some integer, and we introduce sketch tensors \(S_{v_j \to k} \colon \prod_{l \in v_j \to k \cap [d] }[n_l] \times [\tilde{r}_{e_j}]\). By contracting \Cref{eqn: core-determining equation} with \(\bigotimes_{j = 1}^{\mathrm{deg}(k)}S_{v_j \to k}\), one obtains the following sketched linear equation for \(G_k\):
\begin{equation}\label{eqn: TTN sketched linear equation}
    B_{k}(i_k, \beta_{\mathcal{E}(k)}) = \sum_{\alpha_{\mathcal{E}(k)}}\prod_{j = 1}^{\mathrm{deg}(k)}A_{v_j \to k}(\beta_{e_j}, \alpha_{e_j})G_k(i_k, \alpha_{\mathcal{E}(k)}),
\end{equation}
where \(B_{k}\) is the contraction of \(D\) with \(\bigotimes_{j = 1}^{\mathrm{deg}(k)}S_{v_j \to k}\), and each \(A_{v_j \to k}\) is the contraction of \(D_{v_j \to k}\) with \(S_{v_j \to k}\). Specifically, the contraction with each \(S_{v_j \to k}\) is done by contracting the \(i_{v_j \to k \cap [d]}\) multi-index. 
Writing \(A_k = \bigotimes_{j = 1}^{\mathrm{deg(k)}}A_{v_j \to k}\), one sees that \Cref{eqn: TTN sketched linear equation} can be simplified to the linear equation \[\sum_{\alpha_{\mathcal{E}(k)}}A_{k}(\beta_{\mathcal{E}(k)}, \alpha_{\mathcal{E}(k)} )G_k(\alpha_{\mathcal{E}(k)}, i_k) = B_{k}(\beta_{\mathcal{E}(k)}, i_k).\]

The case for \(k = d + 1, \ldots, \tilde{d}\) is similar. We have
\begin{equation*}
    D(i_1, \ldots, i_d) = \sum_{a_{\mathcal{E}(k)}}G_{k}(a_{\mathcal{E}(k)}) \prod_{j = 1}^{\mathrm{deg(k)}}D_{v_j \to k}(\alpha_{e_j},i_{v_j \to k \cap [d]}),
\end{equation*}
for which we apply the contraction with \(\bigotimes_{j = 1}^{\mathrm{deg}(k)}S_{v_j \to k}\). The formula for \(A_{v_j \to k}\) does not change, and the formula for \(B_k\) is
\begin{equation}\label{eqn: TTN sketched linear equation internal}
    B_{k}(\beta_{\mathcal{E}(k)}) = \sum_{\alpha_{\mathcal{E}(k)}}\prod_{j = 1}^{\mathrm{deg}(k)}A_{v_j \to k}(\beta_{e_j}, \alpha_{e_j})G_k(i_k, \alpha_{\mathcal{E}(k)}).
\end{equation}

Thus, for \(k = d + 1, \ldots, \tilde{d}\), the sketched linear equation is
\[\sum_{\alpha_{\mathcal{E}(k)}}A_{k}(\beta_{\mathcal{E}(k)}, \alpha_{\mathcal{E}(k)} )G_k(\alpha_{\mathcal{E}(k)}) = B_{k}(\beta_{\mathcal{E}(k)}).\]

\paragraph{Obtain \(B_k\) by samples}
We shall solve \(G_k\) by the aforementioned linear equation \(A_k G_k = B_k\). One sees that in fact \(B_{k}\) is fully specified. We show how to approximate \(B_k\) by the samples \( \left\{ x^{(j)}:=\left(x_1^{(j)}, \ldots, x_{d}^{(j)} \right) \right\}_{j = 1}^{N}\). 

To perform the sketching with samples as input, we construct each sketch tensor \(S_{v \to k}\) \emph{implicitly} with a continuous sketch function \(s_{v \to k} \colon \R^{|v \to k|} \times [\tilde{r}_{(v, k)}] \to \R\). The function \(s_{v \to k}\) is required to be easy to evaluate, and \(s_{v \to k}\) is related to \(S_{v \to k}\) through the following construction:
\begin{equation*}
    s_{v \to k}(x_{v \to k \cap [d]}, \beta) = \sum_{i_{j} \in [n_j], j \in v \to k \cap [d]} S_{v \to k}(i_{v\to k \cap [d]}, \beta) \prod_{j \in v \to k \cap [d]}\psi_{i_j, j}(x_j).
\end{equation*}

We assume that one is given \(s_{v \to v'}\) over all \(v ,v'\). For \(k = 1,\ldots, d\), one has
\[
\begin{aligned}
    B_k(i_k, \beta_{\mathcal{E}(k)}) = &\int_{\R^d}  p(x_1, \ldots, x_d) \psi_{i_k, k}(x_{k}) \prod_{v \in \mathcal{N}(k) }s_{v \to k}(x_{v \to k \cap [d]}, \beta_{(v, k)})\, d \, w
    \\
    = &\mathbb{E}_{X \sim p}\left[\prod_{v \in \mathcal{N}(k)}s_{v \to k}(X_{v \to k \cap [d]}, \beta_{(v, k)}) \psi_{i_k, k}(X_{k})\right],
\end{aligned}
\]
where the first equality uses the orthonormality of the function basis \(\Vec{\Psi}_j\).

Therefore one can approximate \(B_k\) with samples by
\begin{equation}\label{eqn: tree-based FTN B_k density estimation formula internal node case external}
    B_k(i_k, \beta_{\mathcal{E}(k)}) \approx \frac{1}{N}\sum_{j = 1}^{N}\left(\prod_{v \in \mathcal{N}(k)}s_{v \to k}(x^{(j)}_{v \to k  \cap [d]}, \beta_{(w, k)}) \psi_{i_k, k}(x^{(j)}_{k})\right).
\end{equation}

When \(k = d + 1, \ldots, \tilde{d}\), one likewise obtains the sample-wise approximation
\begin{equation}\label{eqn: tree-based FTN B_k density estimation formula internal node case internal}
    B_k(\beta_{\mathcal{E}(k)}) \approx \frac{1}{N}\sum_{j = 1}^{N}\left(\prod_{v \in \mathcal{N}(k)}s_{v \to k}(x^{(j)}_{v \to k \cap [d]}, \beta_{(v, k)})\right).
\end{equation}

\paragraph{Obtain \(A_k\)}

\Cref{eqn: TTN sketched linear equation} is the linear equation we shall use for \(G_k\), and so one needs to obtain \(B_{k}\) as well as each \(A_{v_j \to k}\). One sees that each \(A_{v_j \to k}\) depends on the gauge degree of freedom for \(D_{v_j \to k}\). Let \(v \in \mathcal{N}(k)\) be one of the \(v_j\). To fully specify \(A_{v \to k}\), we calculate the contraction of \(D\) with \(S_{k \to v} \bigotimes S_{v \to k}\), thereby obtaining
\begin{equation}\label{eqn: equation for Z_k_w}
    Z_{(k, v)}(\beta, \gamma) = \sum_{i_{[d]}}D(i_1, \ldots, i_d)S_{k \to v}(i_{k \to v}, \beta)S_{v \to k}(i_{v \to k}, \gamma).
\end{equation}

One performs the singular value decomposition (SVD) on \(Z_{(k, v)}\), resulting in \(Z_{(k, v)}(\beta, \gamma) = \sum_{\mu}U(\beta, \mu)W(\mu, \gamma)\), and the choice of \(A_{k \to v} = U\) and \(A_{v \to k}= W\) forms a unique and consistent choice of gauge between the pair \((D_{v \to k}, D_{k \to v})\). As obtaining \(U, W\) requires merging the middle factor of the SVD to either the left or right factor, one can use the root information to specify \(U, W\) fully. In particular, if \(v\) is the parent to \(k\), then one takes \(U\) to be the left orthogonal factor in the SVD. If \(k\) is the parent to \(v\), then one takes \(W\) to be the right orthogonal factor in the SVD.

Finally, \(Z_{(k, v)}\) requires a contraction with \(D\), and it can be obtained by samples by 
\begin{equation}\label{eqn: Z_k_v density estimation formula}
    Z_{(k, v)}(\beta, \gamma) \approx \frac{1}{N}\sum_{j = 1}^{N}\left(s_{k \to v}(x^{(j)}_{k \to v\cap [d]}, \beta) s_{v \to k}(x^{(j)}_{v \to k \cap [d]}, \gamma)\right).
\end{equation}

\paragraph{Summary}
We summarize the density estimation algorithm in \Cref{alg:TTN density estimation internal}.

\begin{algorithm}[h]
\caption{Tree-based FTN density estimation.}
\label{alg:TTN density estimation internal}
\begin{algorithmic}[1]
\REQUIRE Sample \(\{x^{(j)}\}_{j = 1}^{N}\).
\REQUIRE Tree \(T = (V, E)\) with \(V = [\tilde{d}]\) and \(V_{\mathrm{ext}} = [d]\). 
\REQUIRE Chosen variable-dependent function basis \(\{\Vec{\Psi}_{j}\}_{j \in [d]}\).
\REQUIRE Collection of sketch tensors \(\{S_{v \to k}, S_{k \to v}\}\) and target internal ranks $\{r_{(k, v)}\}$ for each edge \((k, v) \in E\).
\FOR{each edge \((k, v)\) in \(T\)}
    \STATE Obtain \(Z_{(k,v)}\) by \Cref{eqn: Z_k_v density estimation formula}.
    \STATE Obtain \(A_{k \to v}\) as the left factor of the best rank \(r_{(k, v)}\) factorization of \(Z_{(k,v)}\)
    \STATE Obtain \(A_{v \to k}\) as the right factor of the best rank \(r_{(k, v)}\) factorization of \(Z_{(k,v)}\).
\ENDFOR
\FOR{each node \(k\) in \(T\)}
\IF{node \(k\) contains an external bond}
    \STATE Obtain \(B_{k}\) by \Cref{eqn: tree-based FTN B_k density estimation formula internal node case external}.
\ELSE
    \STATE Obtain \(B_{k}\) by \Cref{eqn: tree-based FTN B_k density estimation formula internal node case internal}.
\ENDIF
    \STATE With \(\mathcal{E}(k) = \{(v_j, k)\}_{j = 1}^{\mathrm{deg}(k)}\), collect each \(A_{v_j \to k}\).
    \STATE Obtain \(G_{k}\) by solving the over-determined linear system \((\bigotimes_{j = 1}^{\mathrm{deg}(k)} A_{v_j \to k})G_{k} = B_{k}\).
\ENDFOR
\end{algorithmic}
\end{algorithm}

\section{Main formulation}\label{sec: main formulation}
This section details the functional hierarchical tensor under a wavelet basis (FHT-W) for 1D and 2D lattice models. \Cref{sec: background on 1D} covers the preliminary information on 1D wavelet transformation. \Cref{sec: model} details the model architecture for FHT-W and how one applies the ansatz to 1D lattice models. \Cref{sec: 2D lattice model main formulation} details how one applies the ansatz to 2D lattice models. \Cref{sec: extensions} discusses extensions of the FHT-W ansatz to other settings.

\subsection{Background on 1D wavelet transformation}\label{sec: background on 1D}
By possibly performing dimension padding, we assume without loss of generality that the dimension \(d\) satisfies \(d = 2^{L}\). A 1D lattice model has a natural ordering, and we denote the variables by \(x = (x_1, \ldots, x_{2^L})\). Typically, after appropriate normalization, one can think of a 1D lattice model as the discretization of a model on the 1D field \(x(a) \colon [0, 1] \to \R\). Thus, one can think of each \(x_j\) for \(j = 1, \ldots, 2^{L}\) as approximating the field \(x\) at length scale \(2^{-L}\).

The wavelet approximation proposes to perform an iterative coarse-graining of the lattice model in \((x_1, \ldots, x_{2^L})\) into orthogonal signals at different length scales. Thus, the input to the transform is data at scale \(2^{-L}\), and the output to the transform is the data at scale \(2^{-l}\) for \(l = L-1, \ldots, 0\). This iterative coarse-graining approach gives rise to the multiscale nature of this work.

We explain the multiresolution wavelet approximation procedure, and more details can be found in \cite{mallat1999wavelet}. For concreteness, we go through the concept with the simple Haar wavelet. In the first step of the coarse-graining process, the wavelet filter transforms the \(x\) variable into two \(2^{L-1}\)-dimensional variables \(y_{L - 1} = (y_{j, L-1})_{j \in [2^{L-1}]}\) and \(c_{L - 1} = (c_{j, L-1})_{j \in [2^{L-1}]}\). which are defined by the following equation
\[
y_{j, L-1} = \frac{x_{2j - 1} + x_{2j}}{\sqrt{2}}, \quad c_{j, L-1} = \frac{x_{2j - 1} - x_{2j}}{\sqrt{2}}.
\]
The variable \(y_{L-1}\) is the scaling coefficient of \(x\) at length scale \(2^{-(L-1)}\), and \(c_{L-1}\) is the detail coefficient of \(x\) at length scale \(2^{-(L-1)}\). One can see that the mapping \(x \to (y_{L-1}, c_{L-1})\) is an invertible and orthogonal transformation.

Repeating the same approach, for \(l = L-1, \ldots, 1\), one performs wavelet transformation on \(y_{l}\) \(2^{l-1}\)-dimensional variables \(y_{l-1}\) and \(c_{l-1}\) given by the following equation
\[
y_{j, l-1} = \frac{y_{2j - 1, l} + y_{2j, l}}{\sqrt{2}}, \quad c_{j, l-1} = \frac{y_{2j - 1, l} - y_{2j, l}}{\sqrt{2}}.
\]

Finally, we note that \(y_1\) is \(2\)-dimensional the procedure ends at \(y_1 \to (y_0, c_0)\). For notational compactness, we denote \(c_{-1}:= y_{0}\).

The wavelet multiresolution approximation transforms the \(x = (x_1, \ldots, x_{2^{L}})\) variable into the variables \(c = (c_{L-1}, \ldots, c_{1}, c_{0}, c_{-1})\), where each \(c_{l}\) is a \(2^{l}\)-dimensional containing the detail of \(x\) at scale \(2^{-l}\), and \(c_{-1} = y_0\) is the coarse-grained approximation of \(x\) at scale \(1\). The mapping \(x \to (c_{k, l})_{l = L-1, \ldots, -1, k = 1, \ldots,  2^{\max(l, 0)}}\) is an invertible and orthogonal transformation.

One can also consider other wavelets. For numerical experiments used in this work, we use the Daubechies D4 wavelet \cite{daubechies1992ten}, which performs the transformation \(y_l\to (y_{l-1}, c_{l-1})\) with a different filter. In the discrete wavelet transform procedure, we use periodic signal extension mode so that any \(2^{l}\)-dimensional variable is transformed to two \(2^{l-1}\)-dimensional variables \cite{mallat1999wavelet}. Extension of this model to pther discrete wavelet transform methods are discussed at the end of this section.

\subsection{Model architecture}\label{sec: model}

\begin{figure}[h]
    \centering
    \includegraphics[width= 0.8
    \linewidth]{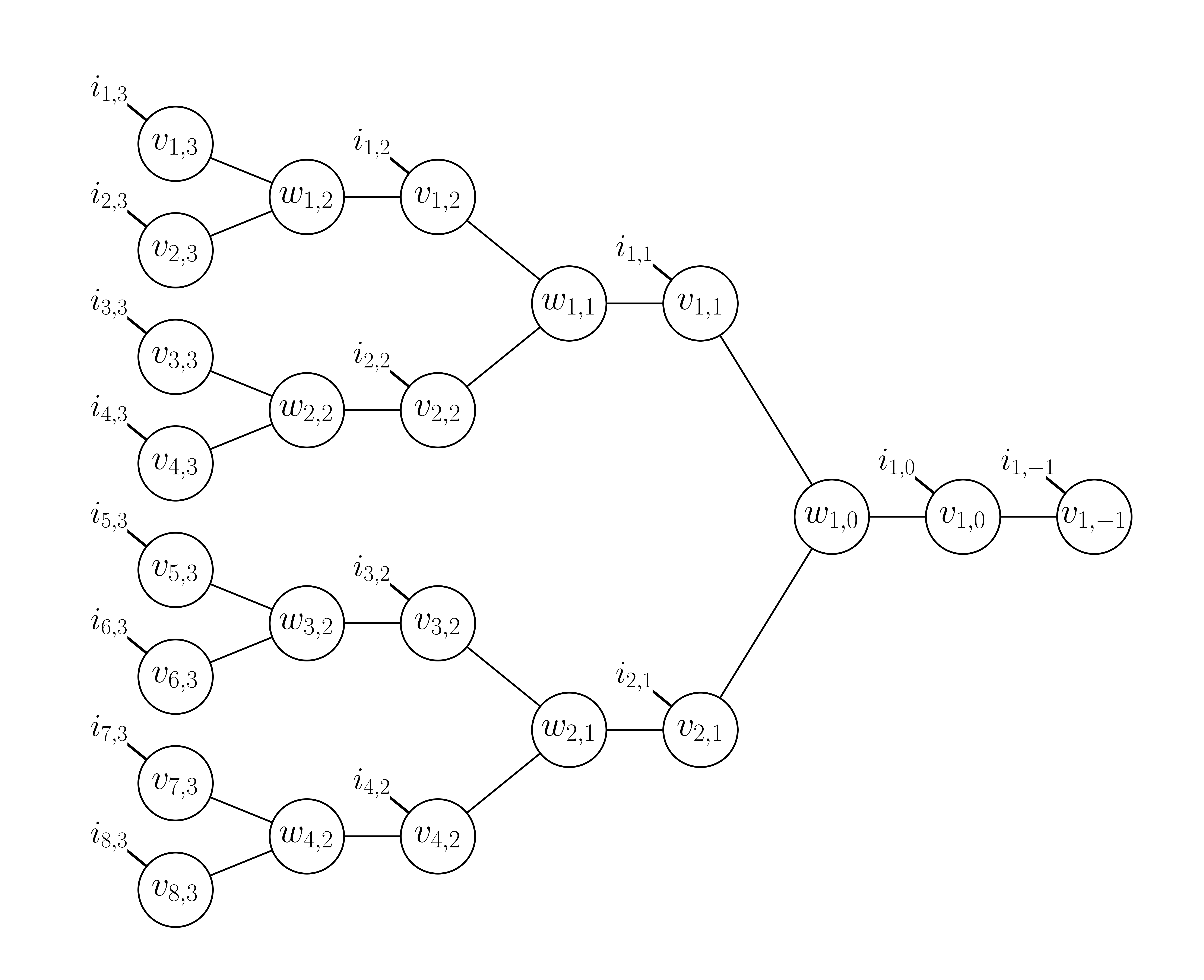}
    \caption{Illustration of the FHT-W ansatz for \(L = 4\). Each \(v_{k, l}\) represents the tensor component \(G_{v_{k, l}}\) at the external node \(v_{k, l}\), and the \(v_{k, l}\) corresponds to the variable \(c_{k, l}\) in the wavelet transformation. Each \(w_{k, l}\) represents the tensor component \(G_{w_{k, l}}\) at the internal node \(w_{k, l}\).}
    \label{fig:wavelet_FHT_L_3}
\end{figure}

For \(d = 2^{L}\), an FHT-W ansatz is a functional tensor network \(f\) taking the variable \(c = (c_{k, l})_{l = L-1, \ldots, -1, k = 1, \ldots,  2^{\max(l, 0)}}\).
We propose FHT-W as a tree-based FTN based on a tree structure \(T = (V, E)\). We illustrate the ansatz in \Cref{fig:wavelet_FHT_L_3}.

We first specify the vertex of \(T\). For the vertex set, we define the external nodes to be \(V_{\mathrm{ext}} = \{v_{k, l}\}_{l \in [L-1] \cup \{0, -1\}, k \in [2^{\max(l, 0)}]}\). The total vertex set is \(V = V_{\mathrm{ext}} \cup \{w_{k, l}\}_{l \in [L-2]\cup \{0\}, k \in [2^l]}\). The vertex \(v_{k, l}\) corresponds to the variable \(c_{k, l}\), and the variables \(w_{k, l}\) is an internal node inserted at each level \(l = 0, \ldots, L-2\).

The edge set of \(T\) is \(E = E_1 \cup E_2 \cup \{(v_{1, 0}, v_{1, -1})\}\), where \(E_1 = \{(w_{k, l}, v_{k, l})\}_{w_{k, l} \in V}\) and \(E_2 = \{(w_{k, l}, v_{2k-1, l+1}), (w_{k, l}, v_{2k, l+1})\}_{w_{k, l} \in V}\). The definition of \(E\) shows that each variable \(c_{2k-1, l+1}\) and \(c_{2k, l+1}\) is connected to the variable \(c_{k, l}\) through an internal node \(w_{k, l}\). This construction for \(T\) ensures that \(c_{k, l}\) is placed closed to \(c_{k', l'}\) when \((k, l)\) is close to \((k', l')\).

The above description fully characterizes the tree \(T = (V, E)\). The FHT-W ansatz \(f\) is a tree-based FTN based on \(T\). The coefficient tensor is \(D\), and \Cref{sec: density estimation} specifies the procedure to obtain \(f\) when one \(f\) is a probability density function and one has access to samples of \(f\).

We explain some of the intuitions behind the design of the FHT-W ansatz. The design of tensor networks for applied sciences is primarily tailored to the problem structure. 
From the analysis of the O-U model and the \(\phi^4\) model in \cite{marchand2023multiscale}, one sees that a large class of lattice models in practice has the property that \(c_{k, l}\) has limited interaction with \(c_{k', l'}\) unless \((k, l)\) is close to \((k', l')\) in both indices. The goal to design \(T\) is to place strongly coupled variables at closer locations on \(T\). Therefore, a binary tree structure in FHT is the optimal shape to place the variables \(c = (c_{k, l})_{k, l}\). To ensure maximal efficiency, we insert internal nodes \(w_{k, l}\), which is done so as to prevent any tensor component in FHT-W from having more than three indices. In practice, a tensor component with four indices is less efficient.

\subsection{Extension to 2D lattice}\label{sec: 2D lattice model main formulation}

We now proceed to the main formulation for 2D lattice models. In this case, we denote the variables by \(x = (x_{(i, j)})_{1 \leq i, j \leq m}\), where \(d = m^2\). Without loss of generality, we assume that \(m = 2^{L}\) so that \(d\) satisfies \(d = 2^{2L}\). In this case, one can think of the 2D lattice model as the discretization of a model on the function field \(x(a) \colon [0, 1]^2 \to \R\). Thus, one can think of each \(x_{(i, j)}\) for \(i,j = 1, \ldots, 2^{L}\) as approximating the field \(x\) at length scale \(2^{-L}\). 

\paragraph{2D wavelet coarse-graining}
In the multidimensional case, we construct the associated filters by separable products of 1D wavelet filters. To apply the iterative wavelet coarse-graining, we perform 1D wavelet filter transformation on \(x = (x_{(i, j)})_{1 \leq i, j \leq m}\) first along the horizontal \(i\) direction followed and then along the vertical \(j\) direction.
As in the 1D case, we consider the periodic signal extension mode, and examples of the 1D wavelet filter include the Haar filter and the D4 filter.  This procedure transforms the \(x\) variable into four \(2^{2L -2 }\)-dimensional variables \(y_{L-1}, c^{lh}_{L -1}, c^{hl}_{L -1}, c^{hh}_{L -1}\), which are respectively obtained by applying (low, low), (low, high), (high, low), (high, high) filters on the \((i, j)\) directions. 

Subsequently, at each \(q = L-1, \ldots, 1\), the 2D wavelet filter transforms \(y_{q}\) into four \(2^{2q-2}\)-dimensional variables \(y_{q-1}\), \(c^{lh}_{q -1}\), \(c^{hl}_{q -1}\), \(c^{hh}_{q -1}\) by applying (low, low), (low, high), (high, low), (high, high) filters on the \((i, j)\) directions. At \(q = 1\), we have \(y_{1}\) transforming into \(y_{0}\), \(c^{lh}_{0 }\), \(c^{hl}_{0 }\), \(c^{hh}_{0}\).

\begin{figure}
    \centering
    \includegraphics[width=\linewidth]{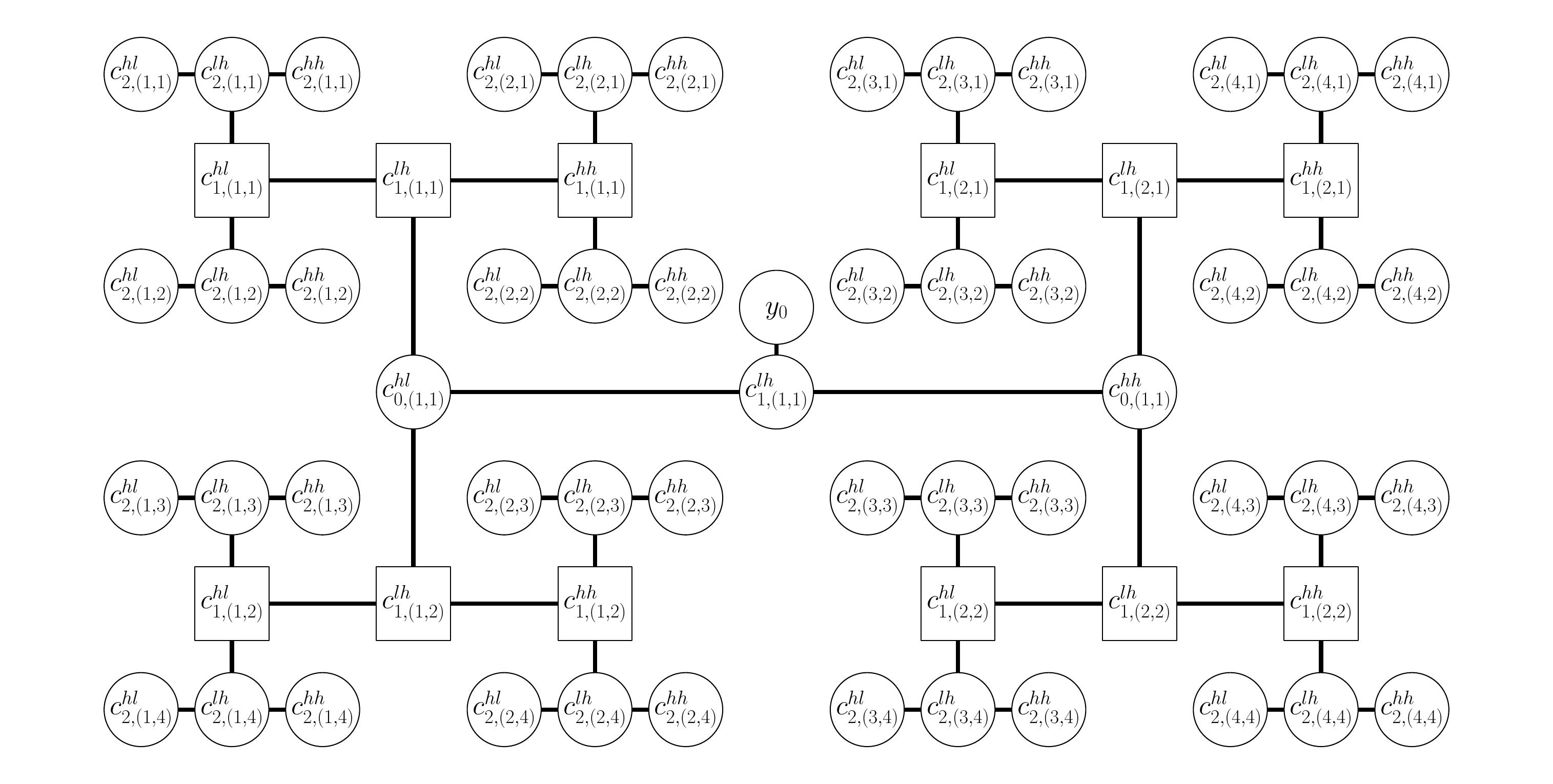}
    \caption{Illustration of the chosen tree structure for 2D wavelet iterative coarsening for a \(d = 8^2\) lattice system. We use circles to illustrate nodes at \(q=2, 0\), and we use squares to illustrate nodes at \(q=1\). Internal nodes on the tree \(T\) are omitted for simplicity.}
    \label{fig:2D_GZ_bipar_alternative}
\end{figure}

\paragraph{Tree architecture}
In this case, we use the tree \(T\) defined in \Cref{sec: model} with \(2L\) levels. We give the 2D wavelet-transformed variables a hierarchical structure \(c = (c_{2L-1}, \ldots, c_{0}, c_{-1})\) with \(c_{l} = (c_{1, l}, \ldots, c_{2^l, l})\), which would allow one to use the FHT-W ansatz in \Cref{sec: model} and the density estimation subroutine in \Cref{alg:TTN density estimation internal}.

We see that each \(c^{lh}_{q}, c^{hl}_{q}, c^{hh}_{q}\) variable is naturally indexed by \((c^{lh}_{q})_{i,j}\), \((c^{hl}_{q})_{i,j}\), \((c^{hh}_{q})_{i,j}\) where \(1 \leq i, j \leq 2^{q}\). At the coarsest level \(l = -1, 0, 1\), we define \(c_{-1} := y_0\), \(c_{0} := c^{lh}_{1, 0}\) and \(c_{1} := (c^{hl}_{0 }, c^{hh}_{0})\). 
On level \(l = 2q, 2q+1\), for each index \(i, j \in [2^q]\), one performs a length \(q\) binary expansion $i = a_{1}\ldots a_{q}$ and $j = b_{1}\ldots b_{q}$. We take \(k\) to be the integer taking the length \(2q\) binary expansion \(k = a_{1} b_{1} \ldots a_{q} b_{q}\). We perform the variable identification in \(c^{lh}_{q}, c^{hl}_{q}, c^{hh}_{q}\) by defining \(c_{k, 2q} := (c^{lh}_{q})_{i,j}, c_{2k-1, 2q+1} = (c^{hl}_{q})_{i,j}\), and \(c_{2k, 2q+1} := (c^{hh}_{q})_{i,j}\). The corresponding tree structure is illustrated in \Cref{fig:2D_GZ_bipar_alternative}.

We give some intuition on how this hierarchical structure is designed. The constructed structure on \(((c^{lh}_{q}, c^{hl}_{q}, c^{hh}_{q})_{q = 0}^{L-1})\) identifies \(c_{2q}\) with \(c^{lh}_{q}\) and identifies \(c_{2q+1}\) with \((c^{hl}_{q}, c^{hh}_{q})\). Moreover, an appropriate flattening is applied so that the aforementioned variable identification remains consistent. Following the variable identification, the illustration in \Cref{fig:2D_GZ_bipar_alternative} shows that \((c^{lh}_{q})_{i,j}\), \((c^{hl}_{q})_{i,j}\) and \((c^{hh}_{q})_{i,j}\) are placed on nodes close to each other on \(T\). Furthermore, this flattening scheme allows \((c^{lh}_{q})_{i,j}\) and \((c^{lh}_{q})_{i',j'}\) to be closer on the binary tree \(T\) if \((i, j)\) is close to \((i', j')\), and the same is true for \(c^{hl}_{q}\) and \(c^{hh}_{q}\). Therefore, the variable identification procedure specified by the flattening procedure gives each variable in \(((c^{lh}_{q}, c^{hl}_{q}, c^{hh}_{q})_{q = 0}^{L-1})\) an appropriate location on \(T\). While the geometric structure of \(T\) cannot reflect all strongly coupled variable pairs, our construction retains the global 2D structure quite well. 

\subsection{Extensions}\label{sec: extensions}
\paragraph{Extension to other wavelet filters}
When one calculates the discrete wavelet transform a variable \(y_{l}\) into two variables \((y_{l-1}, c_{l-1})\), it is in general unclear if the dimension of \((y_{l-1}, c_{l-1})\) matches with that of \(y_{l}\), as the dimension of \((y_{l-1}, c_{l-1})\) might exceed that of \(y_{l}\). For example, if one uses the Daubechies D4 wavelet with the zero padding as the signal extension mode, then \((y_{l-1}, c_{l-1})\) exceeds the dimension of \(y_{l}\) by two due to boundary effects. In this case, one can incorporate these extra variables by inserting more nodes to \(T\) in \Cref{fig:wavelet_FHT_L_3}. 
In addition, for general discrete wavelet transform, one does not require the variable set \(x = (x_1, \ldots, x_d)\) to satisfy \(d = 2^L\). In such cases, one can use an incomplete binary tree by pruning certain branches in \Cref{fig:wavelet_FHT_L_3}. For both cases, one can perform density estimation based on \Cref{sec: density estimation}.

\paragraph{Extensions to the Kolmogorov backward equation}
In the case of FHT, by exclusively using a high-dimensional density estimation subroutine, one can construct the solution operator for the Kolmogorov backward equation (KBE). The approximation of the operator is in an FHT format, and it is termed the operator-valued functional hierarchical tensor \cite{tang2024solvinga}. In particular, when the terminal condition \(f\) in the KBE is in an FHT format, one can use a simple tensor contraction between the FHT ansatz of \(f\) and the FHT ansatz of the approximated solution operator. The same procedure applies in this case. By first applying a 1D or 2D iterative wavelet coarse-graining, one can use the FHT-W ansatz to approximate the solution operator. By constructing the terminal condition \(f\) in an FHT-W format, one can perform an efficient tensor contraction to obtain the approximated PDE solution.

\paragraph{Extension to stochastic optimal control}
The work in \cite{tang2024solvingb} provides a framework for solving a time-dependent Hamilton-Jacobi-Bellman equation through solving a series of stochastic optimal control problems. When one considers optimal control problems in 1D or 2D lattice models, it might be beneficial to use the FHT-W ansatz, and the same numerical treatment in \cite{tang2024solvingb} can be performed in the FHT-W formulation. 

\section{Experiments on numerical rank}\label{sec: numerical rank}
The key motivation for the FHT-W ansatz is that the wavelet transformation gives a better low-rank structure to the target probability distribution \(p\). In this section, we demonstrate through numerical experiments that the probability distribution \(p\) has a substantially lower numerical rank under a wavelet transformation parameterization. For simplicity, we use \(p(x)\) to denote the probability distribution as a function of the original variable \(x\), and we use \(p(c)\) to denote the probability distribution over the wavelet-transformed variable \(c\). \Cref{sec: background on numerical rank} and \Cref{sec: approximation of numerical rank} go through the background and methodology of approximating the numerical rank of \(p(x)\) and \(p(c)\). \Cref{sec: 1D numerical rank} and \Cref{sec: 2D numerical rank} respectively study the numerical rank of \(p(x)\) and \(p(c)\) in 1D and 2D lattice models. \Cref{sec: discussion numerical rank} summarizes the numerical finding.

\subsection{Background on numerical rank}\label{sec: background on numerical rank}

For both 1D and 2D lattice models, we consider a \(d\)-dimensional probability distribution function \(p \colon \R^d \to \R\). When \(p\) is analytic, one can approximate \(p\) under multivariate Chebyshev approximation, and the associated convergence is exponential in the number of basis functions \cite{trefethen2017multivariate}. Similarly, when \(p\) is smooth and periodic, the multivariate Fourier approximation admits a convergence rate faster than any polynomial decay rate \cite{grafakos2008classical}. Thus, up to negligible error, it is mild to assume that \(p\) admits an FTN ansatz under the equation
\begin{equation*}
    p(x_{1}, \ldots, x_{d}) = \left<D, \, \bigotimes_{j=1}^{d} \Vec{\Psi}_{j}(x_j) \right>,
\end{equation*}
where \(D \colon \prod_{j = 1}[n_j] \to \R\) is the coefficient tensor under some appropriate and sufficiently large function basis \(\{\psi_{i,j}\}_{i = 1}^{n_j}\).

The concept of the coefficient tensor \(D\) allows one to consider the rank of the associated function \(p(x)\).
The function \(p(x)\) is said to be of rank \(r\) along the bipartition \([d] = I \cup J\) if there exist functions \(h_{\alpha} \colon \R^{|I|} \to \R, t_{\alpha} \colon \R^{|J|} \to \R\) for  \(\alpha = 1, \ldots, r\) so that the following holds
\[
p(x) = p(x_1, \ldots, x_d) = \sum_{\alpha = 1}^{r}h_{\alpha}(x_I)t_{\alpha}(x_J),
\]
which translates to the equivalent condition that \(D(i_I; i_J)\) is of rank \(r\). 

\paragraph{Numerical rank}
Exact low-rankness is typically unattainable for lattice models. Instead, we use the concept of a numerical rank to study the approximate low-rankness of probability distribution functions on lattice models.
The \emph{numerical rank} of \(p\) is naturally associated with the numerical rank of the unfolding matrix \(D(i_I; i_J)\). We define the numerical rank of \(p\) along \([d] = I \cup J\) as the smallest integer \(r\) so that there exists a rank-\(r\) factorization \(D(i_I; i_J) \approx UV\) with \( \lVert D(i_I; i_J) - UV \rVert \leq \varepsilon \lVert D\rVert\), and we take \(\lVert \cdot \rVert\) to be the matrix \(2\)-norm. 

Likewise, under the wavelet-transformed variable \(c = (c_{k, l})_{k, l}\), the function \(p(c)\) admits an FTN ansatz of the form
\begin{equation*}
    p(c) = \left<F, \, \bigotimes_{k, l} \Vec{\Psi}_{k, l}(c_{k, l}) \right>,
\end{equation*}
where \(F \colon \prod_{k, l}[n_{k, l}] \to \R\) is the coefficient tensor. The numerical rank of \(p(c)\) is defined along an unfolding matrix of \(F\).

\begin{figure}[h]
    \centering
    \includegraphics[width=0.9\linewidth]{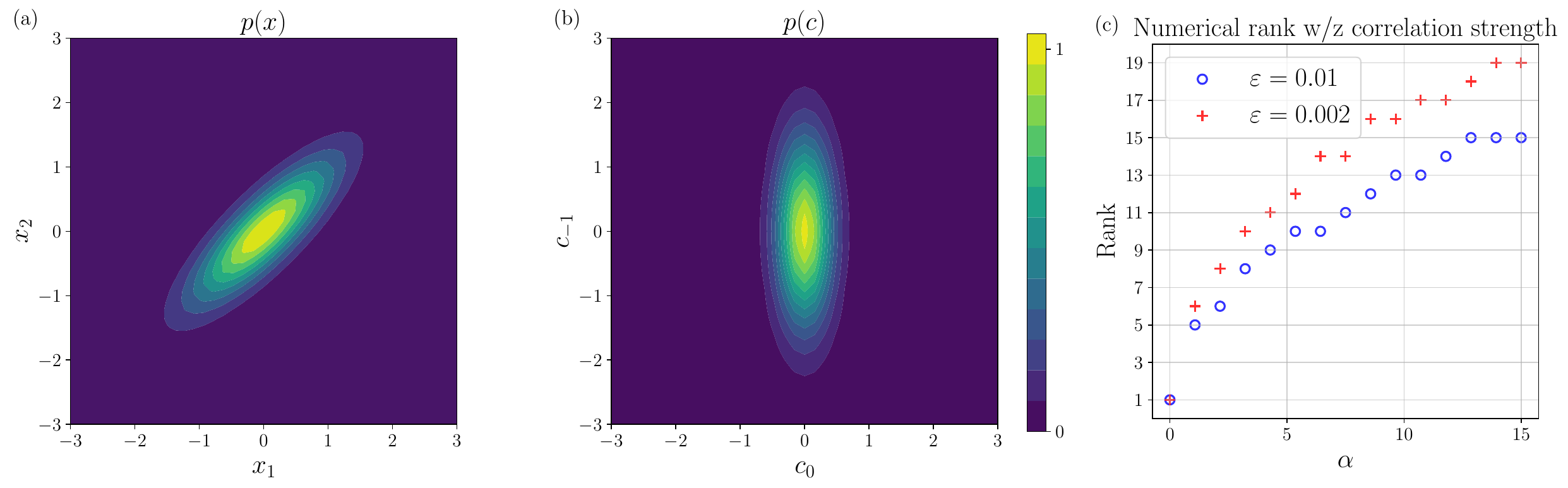}
    \caption{Results of the \(d = 2\) model \(p(x_1, x_2) = \exp\left(-x_1^2/2 - x_2^2/2 - \frac{\alpha}{2}(x_1 - x_2)^2\right)\). \Cref{fig:bivariate_OU_plot_combined}(a)-(b) shows the contour plot of \(p(x)\) and \(p(c)\). \Cref{fig:bivariate_OU_plot_combined}(c) shows that the numerical rank of \(p(x)\) increases with the coupling strength parameter \(\alpha\).}
    \label{fig:bivariate_OU_plot_combined}
\end{figure}

\paragraph{Case study 1: Numerical rank of a \(d = 2\) model} 
We first test the numerical rank on the bivariate model \(p(x_1, x_2) = \exp\left(-\frac{x_1^2}{2} - \frac{x_2^2}{2} - \frac{\alpha}{2}(x_1 - x_2)^2\right)\) with \(\alpha \geq 0\) determining the coupling strength, with larger values of \(\alpha\) leading to stronger coupling. The function \(p\) is the example explored in \Cref{sec: motivating example}. Under the Haar wavelet transform as defined in \Cref{sec: model}, one has one has \((c_{0}, c_{-1}) = (\frac{x_1 - x_2}{\sqrt{2}}, \frac{x_1 + x_2}{\sqrt{2}})\), one sees that \(p(c_{0}, c_{-1}) = \exp(-(\alpha + \frac{1}{2}) c_0^2)\exp( - \frac{c_{-1}^2}{2})\). Thus in this simple case \(p(c)\) is exactly of rank \(1\), which is illustrated by \Cref{fig:bivariate_OU_plot_combined}(b). Therefore, we only need to study the numerical rank of \(p(x)\).

The function \(p(x)\) is effectively supported in the square \(\Omega = [-3, 3]^2\). We let \(\vec{\Psi}_{1}, \vec{\Psi}_{2}\) be the orthogonal Legendre polynomial basis on \(L^2([-3, 3])\). We take the first \(n_1 = n_2 = 60\) Legendre polynomials for each variable, and we obtain the coefficient tensor \(D \in \R^{60 \times 60}\) of \(p\) by Chebyshev interpolation. \Cref{fig:bivariate_OU_plot_combined}(c) shows that the numerical rank of \(D\) increases as \(\alpha\) increases.

\subsection{Approximating the numerical rank}\label{sec: approximation of numerical rank}
For high-dimensional models, it is no longer feasible to perform direct estimation on the numerical rank of the unfolding matrix of \(D, F\), as the rows and columns of the unfolding matrix are typically exponential-sized.
In practice, one can approximate the rank via repeated moment estimations. 
For \(\beta =1, \ldots,  B_1, \gamma = 1, \ldots, B_2\), we let \(s_{I, \beta} \colon \R^{|I|} \to \R \) and \( s_{J, \gamma} \colon \R^{|J|} \to \R\) denote functions respectively taking \(x_{I}\) and \(x_{J}\) as variable. Similar to the procedure in \Cref{sec: density estimation}, one can perform sketching on \(D(i_I; i_J)\) by integrating \(p(x)\) with \(s_{I, \beta}(x_I)s_{J, \gamma}(x_J)\). We define \(Z \in \R^{B_1 \times B_2}\) as the matrix recording all possible integrations of \(p\), i.e.
\[
Z_{I, J}(\beta, \gamma) = \int p(x_1, \ldots, x_d) s_{I, \beta}(x_I)s_{J, \gamma}(x_J) \, dx,
\]
and one can see that \(Z_{I, J}\) has at most rank \(r\) if \(p(x)\) is of rank \(r\). We take the numerical rank of the approximated \(Z_{I, J}\) matrix as an approximation of the numerical rank of \(p(x)\) along \([d] = I \cup J\). By the same procedure, one can approximate the numerical rank of \(F\) by performing repeated moment estimation on \(p(c)\). In other words, for both \(p(x)\) and \(p(c)\), we form \(Z_{I,J}\) and we define the \emph{approximate numerical rank} of \(p(x)\) and \(p(c)\) as the minimal rank for a factorization \(U, V\) to achieve \(\lVert Z_{I, J} - UV \rVert \leq \varepsilon \lVert Z_{I, J}\rVert\). 

To make sure that the numerical rank of \(D, F\) is accurately estimated by the numerical rank of \(Z_{I, J}\), one needs to make sure that the sketching procedure on \(D, F\) is reasonable. For example, one necessary condition is the assumption that the left and right factors of the SVD factorization of \(D, F\) satisfy incoherence conditions. By carrying out the analysis in \cite{tropp2011improved}, under incoherence condition, one can show that taking \(Z_{I, J}\) to be a submatrix of \(D, F\) by random column and row selection ensures that the numerical rank of \(Z_{I, J}\) matches that of \(D, F\). Inspired by this theoretical analysis, in our experiments, we choose \(s_{I, \beta}, s_{J, \gamma}\) so that \(Z_{I, J}\) forms a submatrix of \(D, F\). As the incoherence condition is hard to verify numerically, we choose a large number of functions \(s_{I, \beta}, s_{J, \gamma}\), and we increase the number of functions until the rank estimation saturates. 

Additionally, there is another source of approximation error from obtaining \(Z_{I, J}\) by moment estimation from samples of \(p\). When one performs SVD on \(Z_{I, J}\), the sample estimation error is quite mild, and the estimated spectra follow the Monte-Carlo rate \cite{chen2021spectral}. To address this type of error, we take a large sample size \(N\) and a relatively large rank truncation parameter \(\varepsilon\). In addition, we choose the sample size so that the rank estimation saturates.

\subsection{1D lattice models}\label{sec: 1D numerical rank}
\paragraph{Case study 2: Approximate numerical rank of 1D O–U model}

In this case, we consider the approximate numerical rank of the 1D Ornstein–Uhlenbeck model
\[p(x_1, \ldots, x_d) = \exp\left(-\frac{\alpha}{2}\sum_{i \sim j}(x_i - x_j)^2 - \frac{1}{2}\sum_{i = 1}^{d}x_i^2\right),\] 
where \(\alpha\) controls the correlation strength, and \(i \sim j\) if \(i - j \equiv 1 \mod{d}\). We take \(d = 128\) and \(\alpha = 1000\). The parameter \(\alpha\) is chosen so that neighbor points have strong coupling with \(\mathrm{corr}(X_{i}, X_{i+1}) \approx 1\) while far away points have moderate coupling with \(\mathrm{corr}(X_{i}, X_{i+\ceil{d/2}}) \approx 1/4\). For both \(p(x)\) and \(p(c)\), we form a \(Z_{I, J}\) matrix by moment estimation. For numerical rank of \(Z_{I, J}\), we take \(\varepsilon = 10^{-2}\), and we choose the threshold parameter \(\varepsilon\) so that \(Z_{I, J}\) is of numerical rank \(r=1\) if \(\alpha = 1\). In summary, our experiment reports a numerical rank of \(r = 28\) for \(p(x)\) and a numerical rank of \(r = 7\) for \(p(c)\). We detail the methodology in subsequent paragraphs.

\begin{figure}[h]
    \centering
    \includegraphics[width=0.9\linewidth]{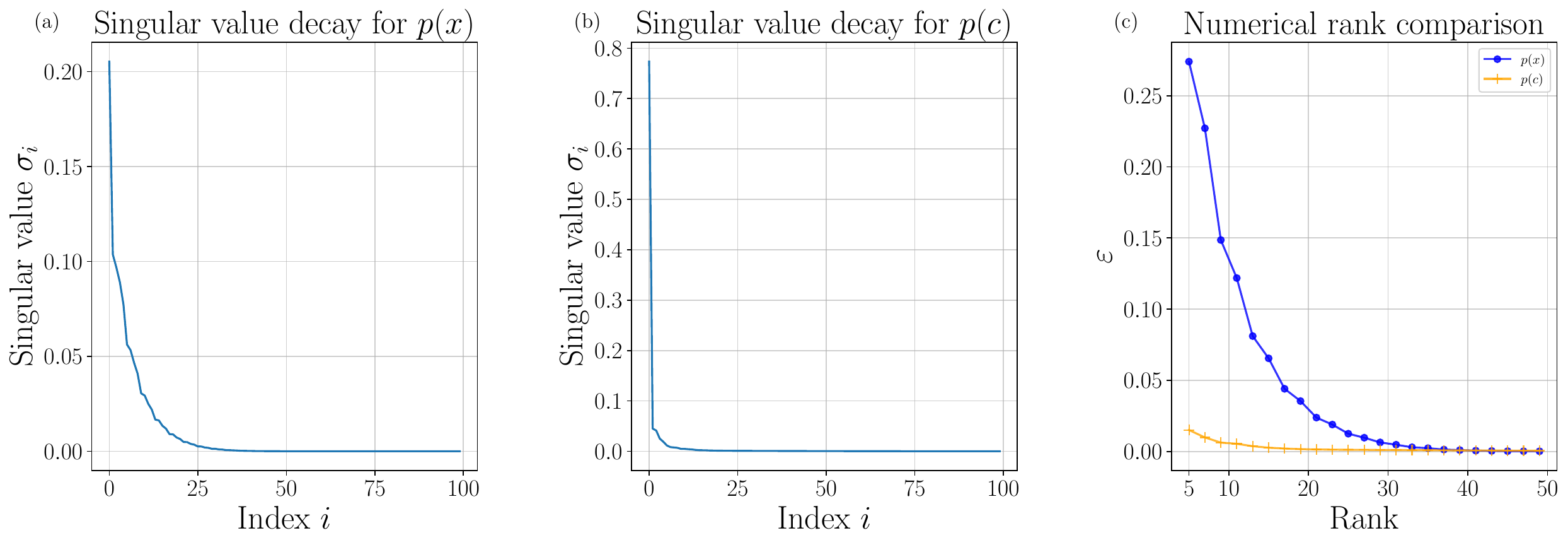}
    \caption{Singular value of \(Z_{I, J}\) for the 1D Ornstein–Uhlenbeck model under \(p(x)\) and \(p(c)\), where \(x\) is the original variable and \(c\) is the variable after wavelet transformation. Both matrices are scaled so that the singular values sum to one, and only the first 100 singular values are plotted. Wavelet transformation leads to more rapid singular value decay in \(Z_{I, J}\).}
    \label{fig:plot_1D_OU}
\end{figure}

We first test the approximate numerical rank of \(p(x)\).
When one uses the FHT and FTT ansatz in modeling the 1D O-U model, one implicitly assumes that the numerical rank of \(p(x)\) is small along the variable bipartition with \(I = \{1, \ldots, d/2\}\) and \(J = \{d/2+1 ,\ldots, d\}\).
By the choice of parameter, the distribution \(X \sim p(x)\) is effectively supported on \(\Omega = [-0.8, 0.8]^d\). 
For \(k = 1, \ldots, d\), we take each \(\vec{\Psi}_{k} = \{\psi_{i, k}\}_{i = 0}^{n}\) to be the collection of orthonormal Legendre polynomial basis functions on \(L^2([-0.8, 0.8])\). We take \(\psi_{i, k}(x_k)\) to be the unique basis function of maximal degree \(i\). The distribution \(X \sim p(x)\) is such that the variable pairs \((X_1, X_d)\) and \((X_{d/2}, X_{d/2+1})\) are strongly coupled. For \(s_{I, \beta}\), we take \(q = 50\) and we construct \(s_{I, \beta}\) to consist of \(\{\psi_{i, k}\}\) for \(i = 0, \ldots, q\) and \(k \in \{1, d/2\}\).  Similarly, we take \(s_{J, \gamma}\) to consist of \(\{\psi_{i, k}\}\) for \(i = 0, \ldots, q\) and \(k \in \{d/2 + 1, d\}\). We note that \(\Psi_{0, k}\) is a constant function for each \(k\), and so there are in total \(B = 2q+1\) non-duplicate functions. The matrix \(Z_{I, J}\) is estimated by \(N = 10^{6}\) samples of \(p\). The singular values of \(Z_{I, J}\) are plotted in \Cref{fig:plot_1D_OU}(a). In this case, the numerical rank of \(Z_{I, J}\) is \(r = 28\). We remark that \(Z_{I, J}\) is a submatrix of \(D\), and the selection of \(s_{I, \beta}\) and \(s_{J, \gamma}\) is chosen such that the rank saturates.

We then test the approximate numerical rank of \(p(c)\). For the wavelet filter, we use the Daubechies D4 filter with periodic signal extension mode. Under the wavelet iterative coarsening, \(p(c)\) takes \(c = (c_{L-1}, \ldots, c_{0}, c_{-1})\) as the variable with \(c_{l} = (c_{1, l}, \ldots, c_{2^l, l})\) for \(l  = 0, \ldots, L-1\). The FHT-W ansatz is motivated by two empirical claims in the renormalization group theory community (c.f. \cite{marchand2023multiscale}): (a) \(p(c)\) is low-rank across scales, and (b) \(p(c)\) is low-rank across physical sites. We test both claims simultaneously. We choose \(I\) to correspond to all variables \(c_{k, l}\) with \(k = 1, \ldots, 2^{l-1}\) and \(l = 1, \ldots, L-1\), and we choose \(J\) to correspond to \(c_{0}, c_{-1}\) and all variables \(c_{k, l}\) with \(k = 2^{l-1} + 1, \ldots, 2^{l}\) and \(l = 1, \ldots, L-1\). This construction makes sure that \(I\) and \(J\) each contains half of the variables for \(c_l\) at each scale \(2^{-l}\). 

As the goal of this work is to justify the use of wavelet transformations, we construct a larger matrix \(Z_{I, J}\) for testing the approximate numerical rank. Similar to the case of \(p(x)\), the distribution of \(C \sim p(c)\) is effectively supported on a hypercube \(\Omega = \prod_{k, l} [-a_{k, l}, a_{k, l}]\), where \(a_{k, l}\) captures the range of the distribution \(C \sim p(c)\) on \(c_{k, l}\). For each \(c_{k, l}\) variable, we let \(\vec{\Psi}_{k, l} = (\psi_{i, (k, l)})_{i =0}^{n}\) be the collection of orthonormal Legendre polynomial basis functions on \(L^2([-a_{k, l}, a_{k, l}])\), where \(\psi_{i, (k, l)}\) corresponds to the basis function of maximal degree \(i\). We take \(q = 30\). 
For each level \(l = L-1, \ldots, 1\), we let \(s_{I, \beta}\) include all \(\psi_{i, (k, l)}\) for \(i = 0, \ldots, q\) and \(k \in \{1, 2^{l-1}\}\), and we let \(s_{J, \gamma}\) include all \(\psi_{i, (k, l)}\) for \(i = 0, \ldots, q\) and \(k \in \{2^{l-1} + 1, 2^{l}\}\). This choice corresponds to taking the boundary variables in \(c_{l}\) at each level \(l\). The singular values of \(Z_{I, J}\) are plotted in \Cref{fig:plot_1D_OU}(b). In this case, the numerical rank of \(Z_{I, J}\) is \(r = 7\).

In both \(p(x)\) and \(p(c)\), we remark that \(Z_{I, J}\) is constructed by forming a submatrix of the unfolding matrix of \(D\) and \(F\). The choice of \(Z_{I, J}\) for \(p(x)\) is less extensive and might lead to underreporting of the approximate numerical rank of \(D\). However, one can see that the approximate numerical rank of \(p(x)\) is more than double that of \(p(c)\). If one models \(p(x)\) with the FHT ansatz and \(p(c)\) with the FHT-W ansatz, the difference in the rank \(r\) leads to an order of magnitude difference in the parameter size due to the \(O(r^3)\) scaling for both models. Moreover, we see from \Cref{fig:plot_1D_OU} that the singular value decay is more rapid for \(p(c)\), which suggests that the sketching algorithm for the FHT-W ansatz is more robust to a small target internal rank parameter \(\{r_{e}\}_{e \in E}\). 

\paragraph{Case study 3: Approximate numerical rank of 1D Ginzburg-Landau model}

In this case, we consider the approximate numerical rank of the 1D Ginzburg–Landau model
\[p(x_1, \ldots, x_d) = \exp\left(-\frac{\alpha}{2}\sum_{i \sim j}(x_i - x_j)^2 - \frac{\lambda}{2}\sum_{i = 1}^{d}(1-x_i^2)^2\right),\] 
where \(\alpha\) controls the correlation strength, \(\lambda\) controls the strength of the double-well term, and \(i \sim j\) if \(i - j \equiv 1 \mod{d}\). We test an interesting case with a strong site-wise interaction and a non-trivial double-well phenomenon, and we choose \(d = 128\), \(\alpha = 250\) and \(\lambda = 5\). We compare the numerical rank between \(p(x)\) and \(p(c)\). 

\begin{figure}[h]
    \centering
    \includegraphics[width=0.9\linewidth]{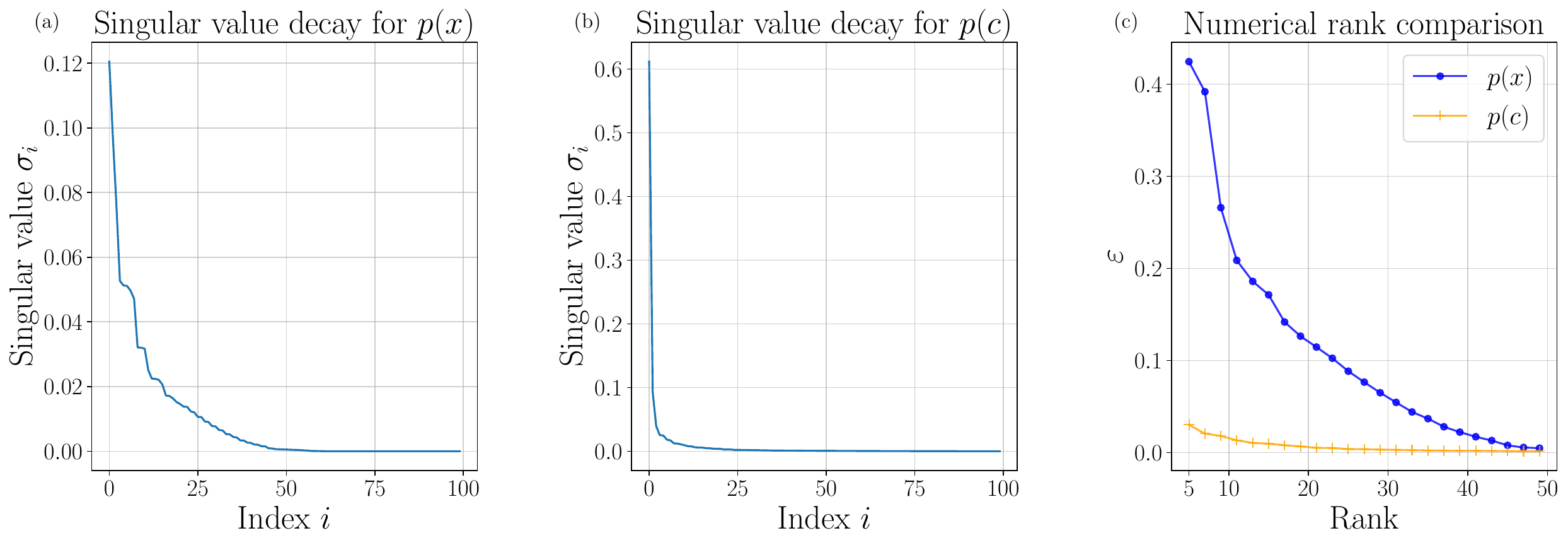}
    \caption{Singular value of \(Z_{I, J}\) for the 1D Ginzburg-Landau model under \(p(x)\) and \(p(c)\). Both matrices are scaled so that the singular values sum to one, and only the first 100 singular values are plotted. Wavelet transformation leads to more rapid singular value decay in \(Z_{I, J}\).}
    \label{fig:plot_1D_GL}
\end{figure}

In this case, we can similarly identify a hypercube support for \(X \sim p(x)\) and \(C \sim p(c)\), and subsequently we can give \(p(x)\) and \(p(c)\) an FTN structure by a Legendre polynomial basis. Thus, by changing the definition of each \(\psi_{i, k}\) and \(\psi_{i, (k, l)}\) to be based on the data support in the G-L model, we can use the same choice of \(s_{I, \beta}\) and \(s_{J, \gamma}\) as in the O-U model.
In this case, we also obtain \(Z_{I, J}\) with moment estimation, and we use \(N = 5 \times 10^5\) samples generated from Markov-chain Monte-Carlo (MCMC) simulations based on \(p\). The result is shown in \Cref{fig:plot_1D_GL}. With \(\varepsilon = 0.01\), the numerical rank for \(Z_{I,J}\) is \(r = 45\) for \(p(x)\), and the numerical rank for \(Z_{I,J}\) is \(r = 14\) for \(p(c)\).

\subsection{2D lattice models}\label{sec: 2D numerical rank}
\paragraph{Case study 4: Approximate numerical rank of 2D O-U model}

In this case, we consider the approximate numerical rank of the 2D Ornstein–Uhlenbeck model,
\[
\begin{aligned}
    &p(x_{(1, 1)}, \ldots, x_{(m, m)}) 
    \\
    =&\exp\left(-\frac{\alpha_1}{2}\sum_{j}\sum_{i \sim i'}(x_{(i, j)} - x_{(i', j)})^2 -\frac{\alpha_2}{2}\sum_{i}\sum_{j \sim j'}(x_{(i, j)} - x_{(i, j')})^2
- \frac{1}{2}\sum_{i,j = 1}^{m}x_{(i,j)}^2\right),
\end{aligned}
\] 
where \(\alpha_1, \alpha_2\) control the correlation strength, and \(i \sim i'\) if \(i - i' \equiv 1 \mod{m}\). We take \(m = 8\) with \(d = m^2 = 64\). When both \(\alpha_1\) and \(\alpha_2\) are large, the rank of \(p(x)\) is not significantly larger than 1D O-U models. Thus, we consider an interesting case where \(\alpha_1 = 200\) and \(\alpha_2 = 10\), which corresponds to the case where the interaction at the first spatial index is considerably stronger than the second spatial index. In summary, the approximate numerical rank of \(p(x)\) is \(r = 73\) and the approximate numerical rank for \(p(c)\) is \(r = 17\).

We first choose the unfolding matrix for \(p(x)\) and \(p(c)\). For \(p(x)\), we take \(I\) to be the collection of \(x_{(i, j)}\) with \(i \in \{1,\ldots, m/2\}\) and \(j \in [m]\), and we take \(J\) to be the collection of \(x_{(i, j)}\) with \(i \in \{m/2+1,\ldots, m\}\) and \(j \in [m]\). For \(p(c)\), we perform a 2D wavelet transformation to get the \(c\) variable, and we take the separable Daubechies D4 filter with periodic signal extension mode. 
The recursive binary expansion detailed in \Cref{sec: 2D lattice model main formulation} leads to the variable structure that is similar to the 1D case, and one has \(c = (c_{2L-1}, \ldots, c_{0}, c_{-1})\) as the variable with \(c_{l} = (c_{1, l}, \ldots, c_{2^l, l})\) for \(l = 1, \ldots, 2L-1\). We choose \(I\) to correspond to all variables \(c_{k, l}\) with \(k = 1, \ldots, 2^{l-1}\) and \(l = 1, \ldots, 2L-1\), and we choose \(J\) to correspond to \(c_{0}, c_{-1}\) and all variables \(c_{k, l}\) with \(k = 2^{l-1} + 1, \ldots, 2^{l}\) and \(l = 1, \ldots, 2L-1\).

Similar to the case of the 1D O-U model, both \(p(x)\) and \(p(c)\) are respectively supported on a hypercube, and so one can use a Legendre polynomial basis function on each variable. Thus, to form \(Z_{I, J}\), we choose \(s_{I, \beta}\) and \(s_{J, \gamma}\) to be univariate Legendre basis functions. For \(p(x)\), we take Legendre polynomial basis functions up to degree \(q = 50\). We let \(s_{I, \beta}\) include the basis functions for \(x_{(i,j)}\) for \(i \in \{1, m/2\}\) and \(j \in [m]\). We let \(s_{J, \gamma}\) include the basis functions for \(x_{(i,j)}\) for \(i \in \{m/2 + 1, m\}\) and \(j \in [m]\). For \(p(c)\), picking the variables of interest is difficult. Therefore, we construct \(s_{I, \beta}\) and \(s_{J, \gamma}\) to include all variables in \(I\) and \(J\) for Legendre polynomial basis functions up to degree \(q = 20\).

\begin{figure}[h]
    \centering
    \includegraphics[width=0.9\linewidth]{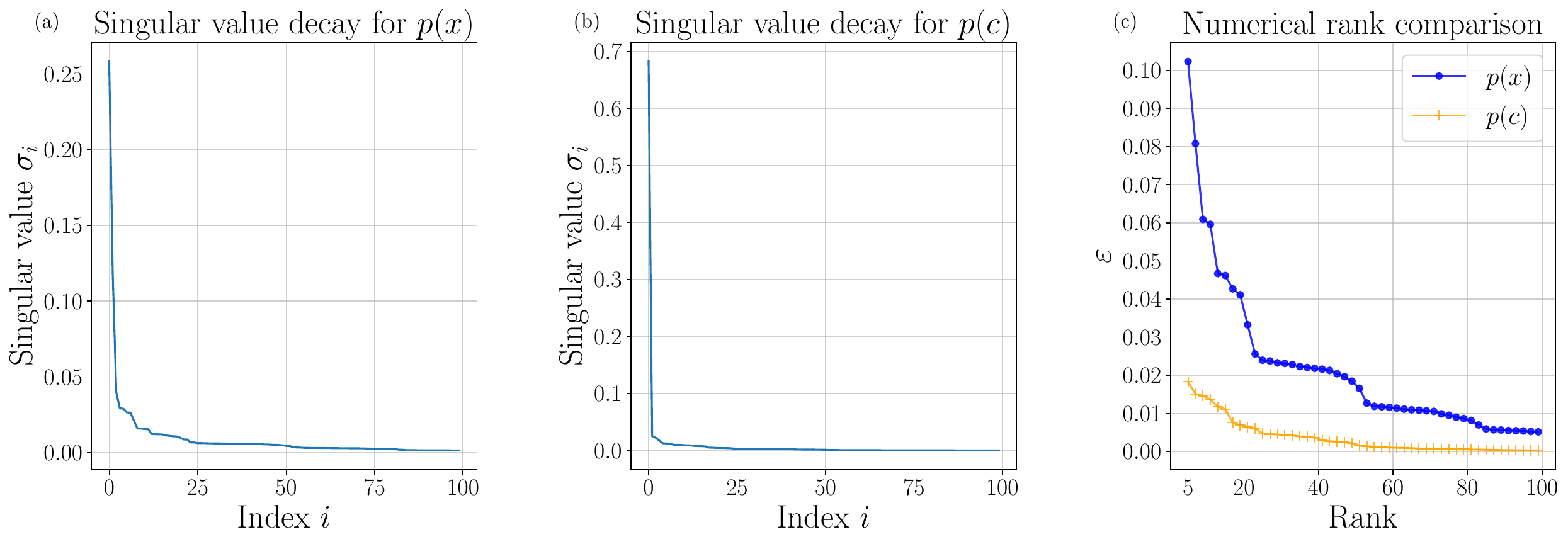}
    \caption{Singular value of \(Z_{I, J}\) for the 2D Ornstein–Uhlenbeck model under \(p(x)\) and \(p(c)\). Both matrices are scaled so that the singular values sum to one, and only the first 100 singular values are plotted. Wavelet transformation leads to more rapid singular value decay in \(Z_{I, J}\).}
    \label{fig:plot_2D_OU}
    \centering
    \includegraphics[width=0.9\linewidth]{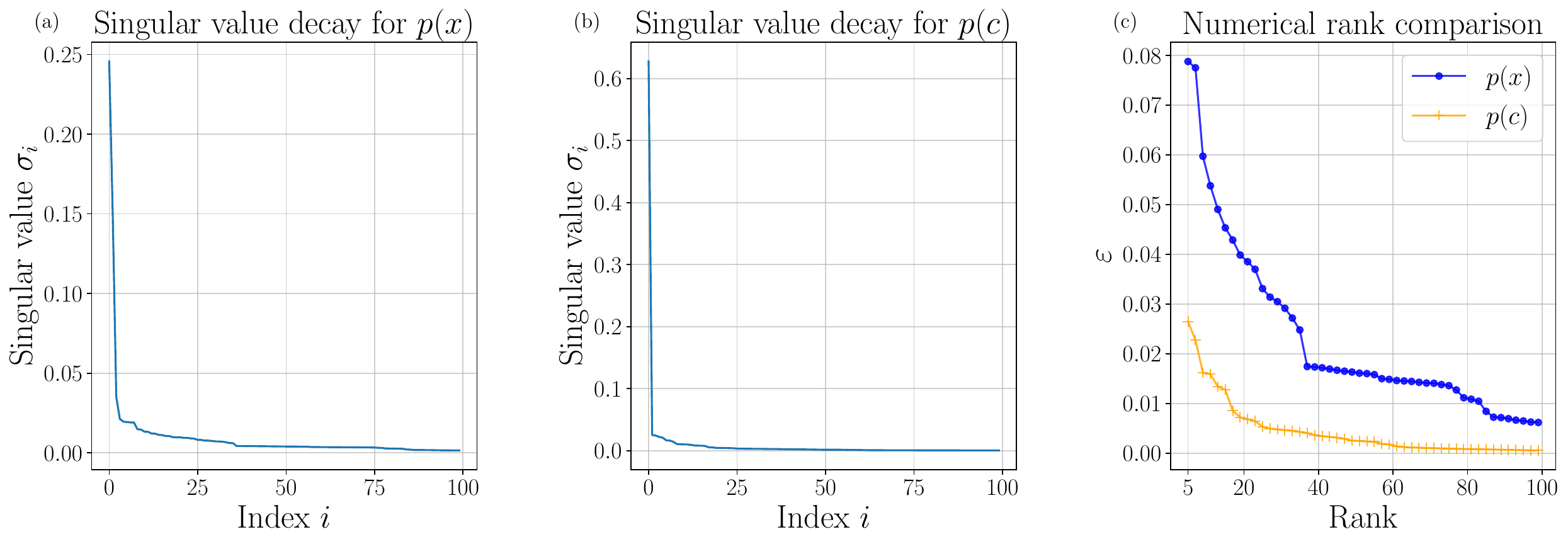}
    \caption{Singular value of \(Z_{I, J}\) for the 2D Ginzburg-Landau model under \(p(x)\) and \(p(c)\). Both matrices are scaled so that the singular values sum to one, and only the first 100 singular values are plotted. Wavelet transformation leads to more rapid singular value decay in \(Z_{I, J}\).}
    \label{fig:plot_2D_GL}
\end{figure}

We use \(N = 5 \times 10^5\) samples to estimate \(Z_{I, J}\) and we take \(\varepsilon = 0.01\). The approximate numerical rank of \(p(x)\) is \(r = 73\) and the approximate numerical rank for \(p(c)\) is \(r = 17\). The results are summarized in \Cref{fig:plot_2D_OU}.

\paragraph{Case study 5: Approximate numerical rank of 2D Ginzburg-Landau model}

In this case, we consider the approximate numerical rank of the 2D Ginzburg–Landau model,
\[
\begin{aligned}
    p(x_{(1, 1)}, \ldots, x_{(m, m)}) 
    =&\exp\Big(-\frac{\alpha_1}{2}\sum_{j}\sum_{i \sim i'}(x_{(i, j)}   - x_{(i', j)})^2  \\ & \,\,\,\,\,\,\,\,\,\,\,\,\,\,\,-\frac{\alpha_2}{2}\sum_{i}\sum_{j \sim j'}(x_{(i, j)} - x_{(i, j')})^2
- \frac{\lambda}{2}\sum_{i,j = 1}^{m}\left(1-x_{(i,j)}^2\right)^2\Big),
\end{aligned}
\] 
where \(\alpha_1, \alpha_2\) control the correlation strength, \(\lambda\) controls the strength of the double-well term, and \(i \sim i'\) if \(i - i'\equiv 1 \mod{m}\). We take \(m = 8\) with \(d = m^2 = 64\). For the model parameter, we take \(\alpha_1 = 20\), \(\alpha_2 = 0.6\) and \(\lambda = 1\). By adjusting the support for each Legendre polynomial basis, we can use the same choice of \(s_{I, \beta}\) and \(s_{J, \gamma}\) as in the 2D O-U model.
We obtain \(Z_{I, J}\) with moment estimation, and we use \(N = 5 \times 10^5\) samples generated from Markov-chain Monte-Carlo (MCMC) simulations based on \(p\). 
The result is shown in \Cref{fig:plot_2D_GL}.
With \(\varepsilon = 0.01\), the numerical rank for \(Z_{I,J}\) is \(r = 84\) for \(p(x)\), and the numerical rank for \(Z_{I,J}\) is \(r = 17\) for \(p(c)\).

\subsection{Discussion of results}\label{sec: discussion numerical rank}
From the experiments carried out in 1D and 2D systems, we show that \(p(c)\) has a more favorable numerical rank than \(p(x)\) across a large collection of lattice models with strong interactions. We make several remarks based on the findings.

The experiment finding indicates when it is appropriate to use the FHT-W ansatz.
In our case studies, the problem parameters are chosen to be numerically challenging instances so that one can draw empirical observations among relatively difficult probability distributions among each category. 
Overall, the case studies lead to three empirical observations for \(p(x)\): (a) Numerical rank of \(Z_{I,J}\) tend to increase with coupling strength, (b) Numerical rank of \(Z_{I, J}\) is typically larger for 2D lattice models than for 1D systems, and (c) Numerical rank tend to increase under the presence of double-well terms. Thus, when one encounters a strongly coupled 2D Ginzburg-Landau model or a \(\phi^4\) model \cite{zinn2021quantum}, it might be preferable to use the FHT-W ansatz over the FHT ansatz defined in \cite{tang2024solving}. 

Moreover, our case study shows that the FHT-W density estimation is more efficient than the FHT ansatz in modeling lattice models. In the 2D G-L model, the rank for \(p(c)\) is smaller by a factor of five, which leads to a difference of two orders of magnitude in the model complexity. For the density estimation subroutine in \Cref{alg:TTN density estimation internal}, the main computational bottleneck lies in performing the moment estimation subroutine. When one has \(N\) samples, obtaining \(B_k\) by density estimation has a cost of \(O(N \tilde{r}^3)\), where \(\tilde{r} := \max_{e} \tilde{r}_{e}\), and we take \(\tilde{r}\) to be larger than the numerical rank to ensure the sketched linear system is well-conditioned. Thus, as the numerical rank is as large as \(r \approx 80\) in the 2D G-L model, we can see that it is numerically challenging to model \(p(x)\) with the FHT density estimation subroutine. In contrast, we see that \(p(c)\) is consistently low-rank with \(r < 20\) in all of the experiments, which shows that the FHT-W ansatz admits a more efficient density estimation subroutine for the studied lattice models.

\section{Numerical experiment}\label{sec: numerical experiment}

This section demonstrates the numerical performance of the proposed density estimation algorithm.
In \Cref{sec: 1D lattice model result}, we test the FHT-W ansatz in density estimation of 1D lattice models. In \Cref{sec: 2D lattice model result}, we test the FHT-W ansatz in density estimation of 2D lattice models.

\subsection{1D lattice models}\label{sec: 1D lattice model result}

\begin{figure}[h]
  \centering
  \includegraphics[width = 0.8\linewidth]{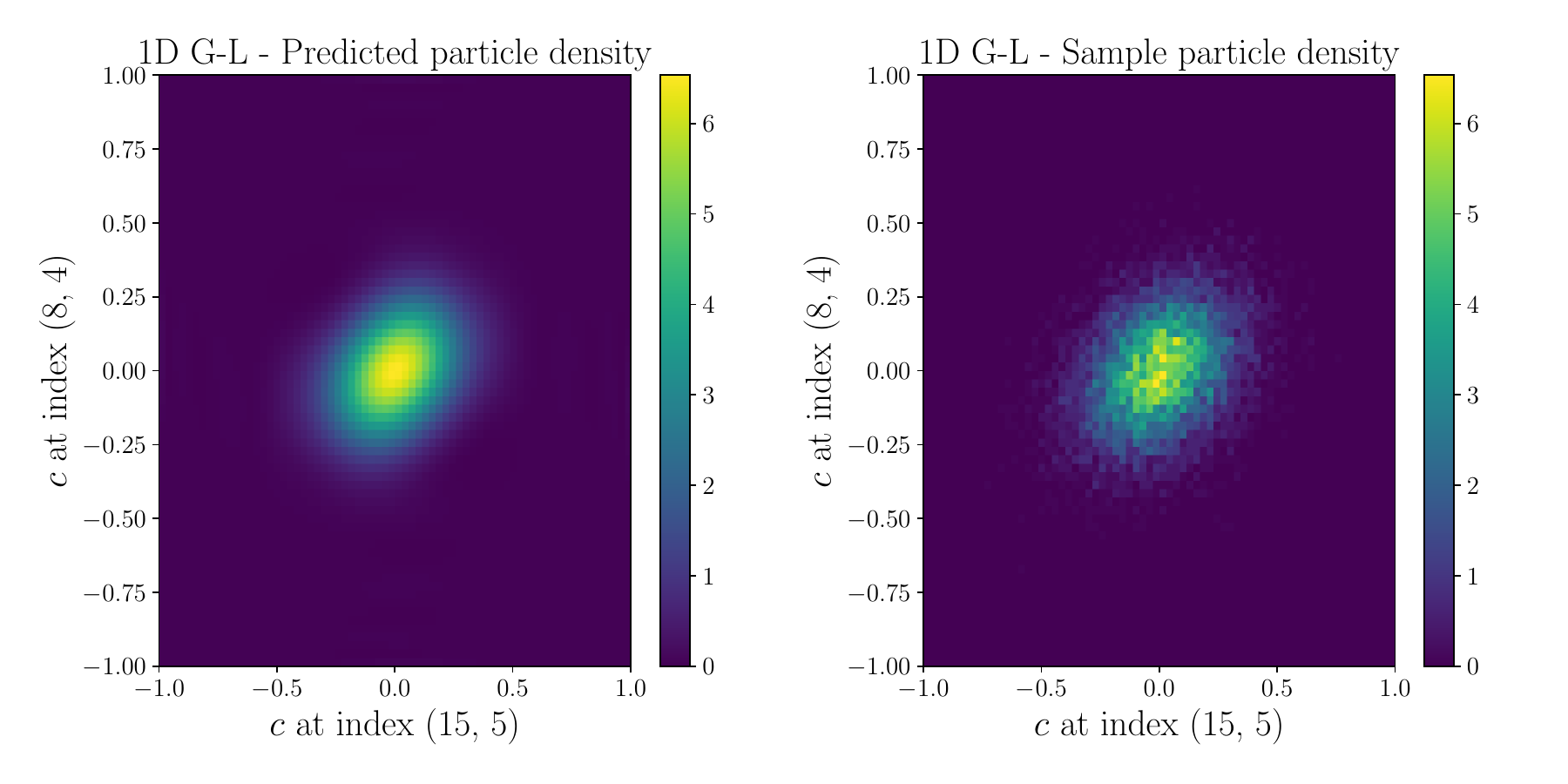}
  \includegraphics[width = 0.8\linewidth]{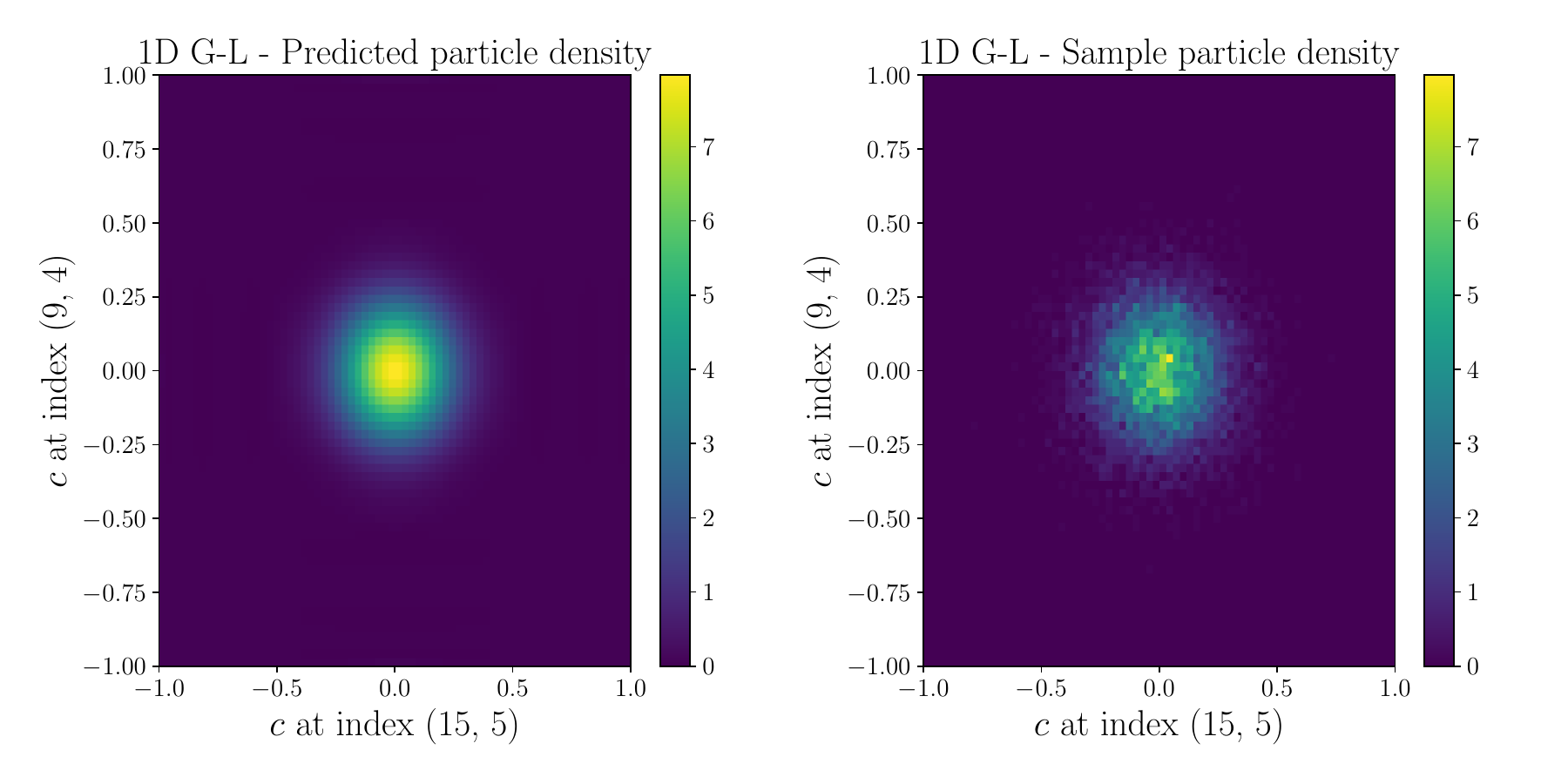}
  \caption{1D Ginzburg-Landau model. Plots of the marginal distribution of \(C \sim p(c)\) at \((c_{15, 5}, c_{8, 4})\) and \((c_{15, 5}, c_{9, 4})\). For illustration purposes, a scaling is performed so that the marginal distribution lies in \([-1,1]^2\).}
  \label{Fig: 1D marginal}
\end{figure}

\begin{figure}[h]
  \centering
  \includegraphics[width = \linewidth]{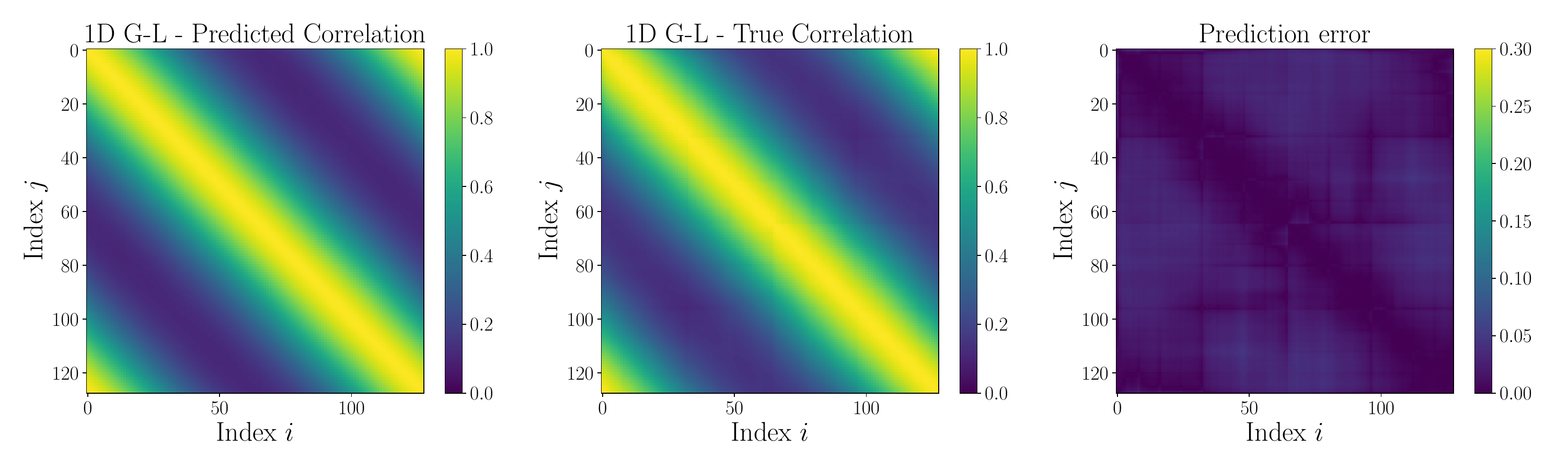}
  \caption{1D Ginzburg-Landau model. Plot of the correlation matrix predicted by the FHT-W ansatz.}
  \label{Fig: 1D covariance}
  \centering
  \includegraphics[width = \linewidth]{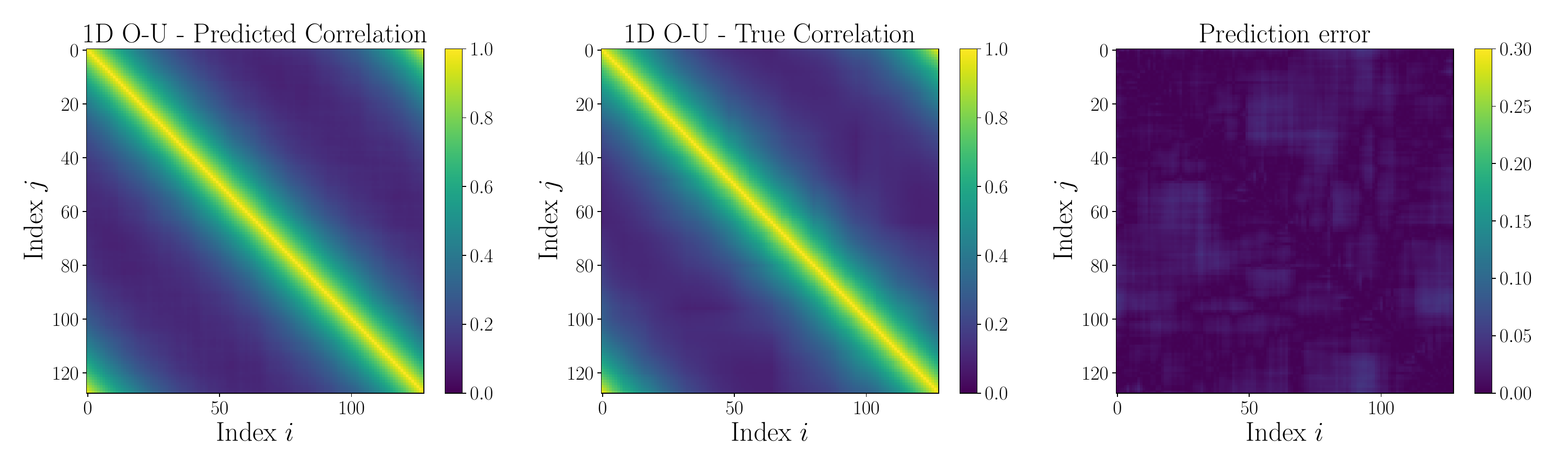}
  \caption{1D Ornstein–Uhlenbeck model. Plot of the correlation matrix predicted by the FHT-W ansatz.}
  \label{Fig: 1D covariance O-U}
\end{figure}

This subsection tests \Cref{alg:TTN density estimation internal} on the density estimation of 1D lattice models. We first consider a strongly coupled 1D Ginzburg-Landau model. The distribution is of the form
\[p(x_1, \ldots, x_d) = \exp\left(-\frac{\alpha}{2}\sum_{i \sim j}(x_i - x_j)^2 - \frac{\lambda}{2}\sum_{i = 1}^{d}(1-x_i^2)^2\right),\] 
where \(i \sim j\) if \(i - j \equiv 1 \mod{d}\). We take the parameter as specified in \Cref{sec: 1D numerical rank} so that \(d = 128, \alpha = 250\) and \(\lambda = 5\). We perform \Cref{alg:TTN density estimation internal} on \(N = 12000\) samples of \(p\) after wavelet transformation. We take the maximal internal dimension to be \(r = 12\). On each variable, we use a Legendre polynomial basis with a maximal degree of \(q = 25\). \Cref{Fig: 1D marginal} shows that the FHT-W ansatz accurately captures the marginal distribution of \(C \sim p(c)\) at \((c_{15, 5}, c_{8, 4})\) and \((c_{15, 5}, c_{9, 4})\). It is noteworthy from the 2-marginal distribution that \((c_{15, 5}, c_{9, 4})\) is already effectively decoupled. Even for the mildly correlated 2-marginal on \((c_{15, 5}, c_{8, 4})\), we see that the distribution is unimodal. The relatively simple behavior of \(p(c)\) at finer scales is one key reason underlying the small numerical rank observed in \Cref{sec: 1D numerical rank}.

To give a challenging numerical test of the FHT-W ansatz, we use the FHT-W ansatz to predict the correlation matrix,
\[
M(i, j) =\mathrm{Corr}_{X \sim p(x)}\left(X_i, X_j\right),
\]
which one can see can be carried out by performing \(O(d^2)\) observable estimation tasks from the FHT-W ansatz obtained from \Cref{alg:TTN density estimation internal}. \Cref{Fig: 1D covariance} shows that the predicted correlation closely matches the ground truth obtained from sample estimation.

We repeat the same numerical experiment on the 1D Ornstein–Uhlenbeck model, where the model parameter likewise follows from that of \Cref{sec: 1D numerical rank}.
The probability density function is
\[p(x_1, \ldots, x_d) = \exp\left(-\frac{\alpha}{2}\sum_{i \sim j}(x_i - x_j)^2 - \frac{1}{2}\sum_{i = 1}^{d}x_i^2\right),\] 
where \(d = 128\) and \(\alpha = 1000\).
We perform \Cref{alg:TTN density estimation internal} on \(N = 12000\) samples of \(p\). In \Cref{Fig: 1D covariance O-U}, we see that the FHT-W ansatz can likewise accurately model the global correlation structure of the 1D O-U model.

\begin{figure}[h]
  \centering
  \includegraphics[width = 0.8\linewidth]{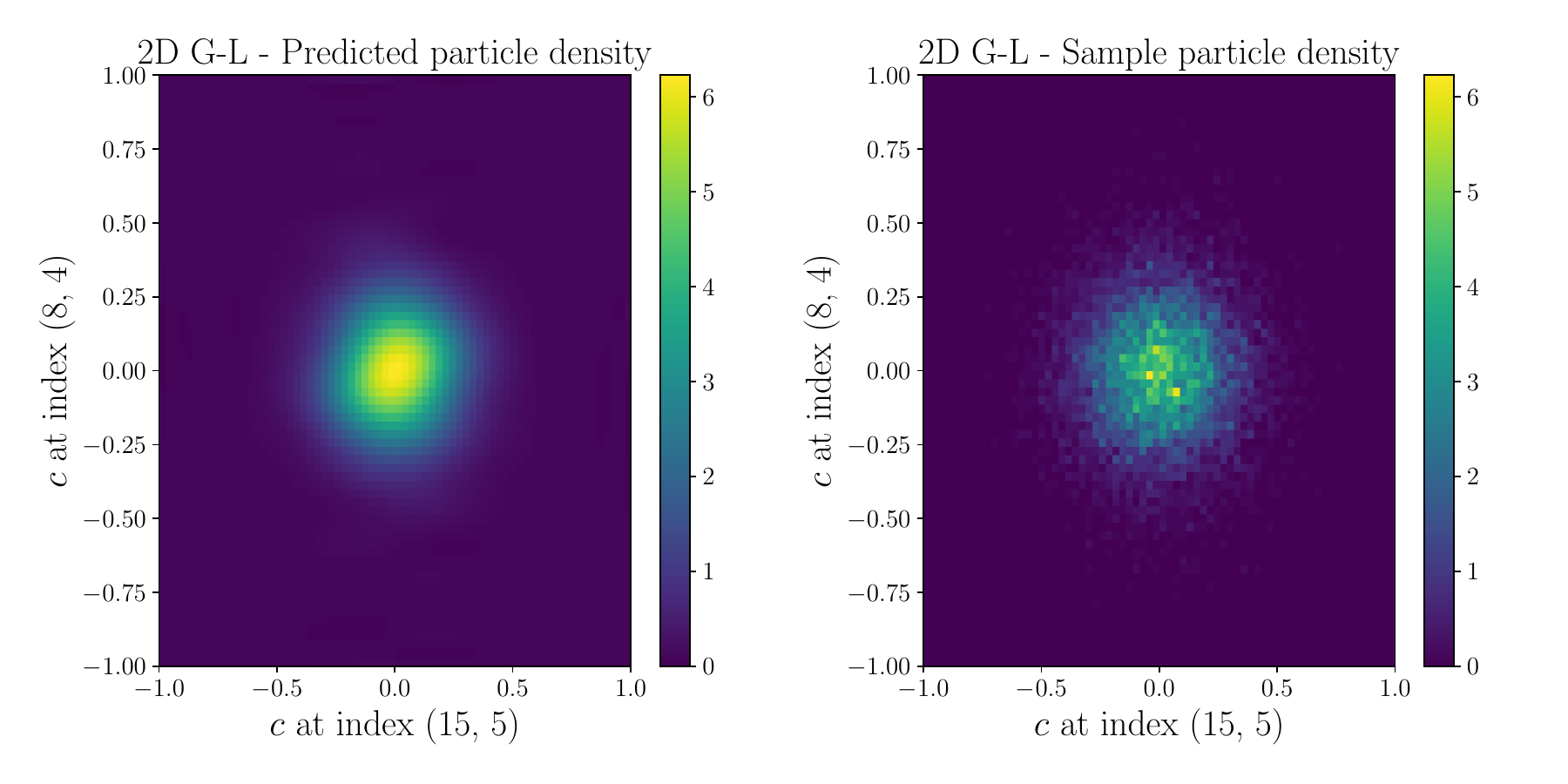}
  \includegraphics[width = 0.8\linewidth]{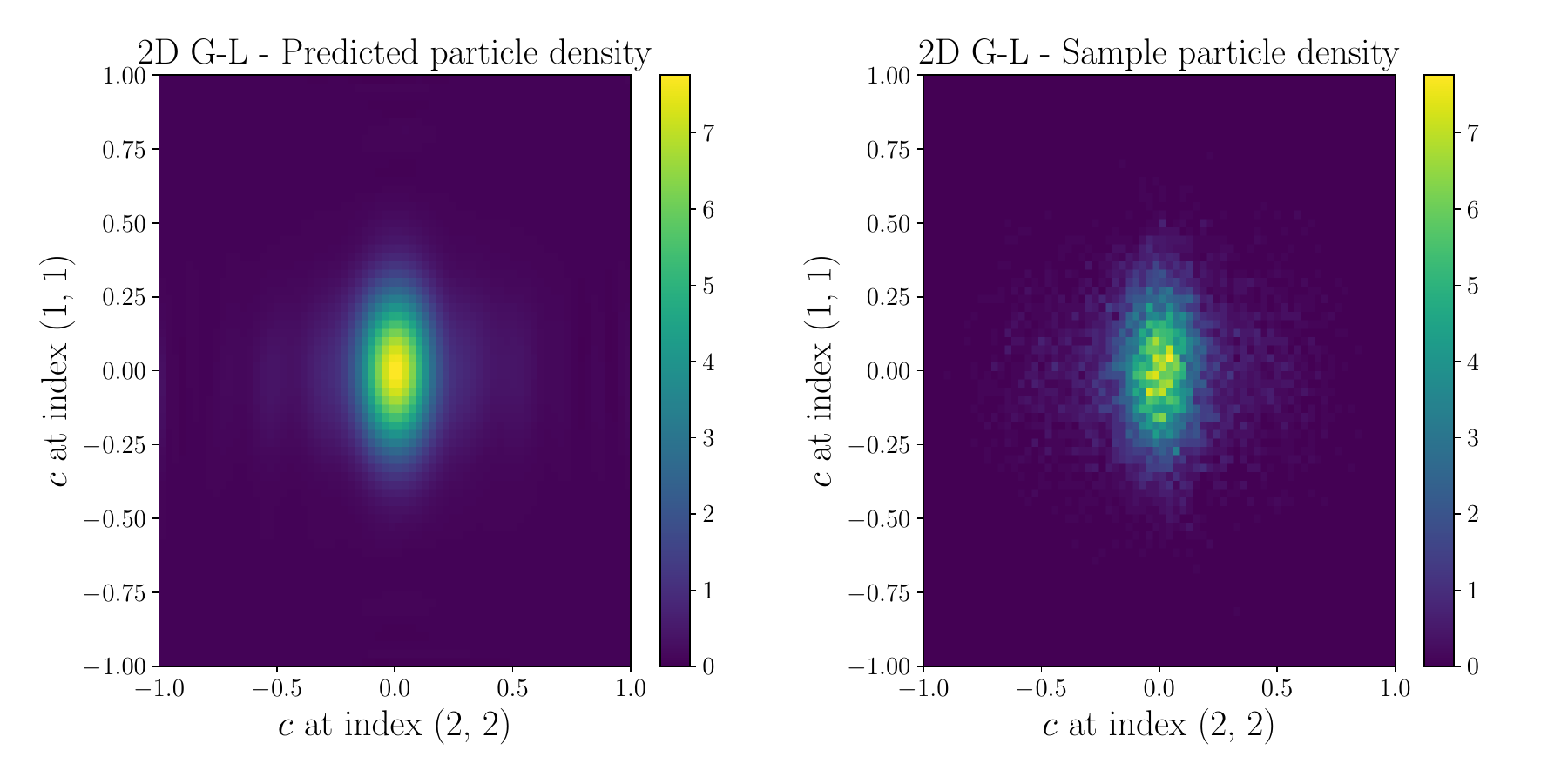}
  \includegraphics[width = 0.8\linewidth]{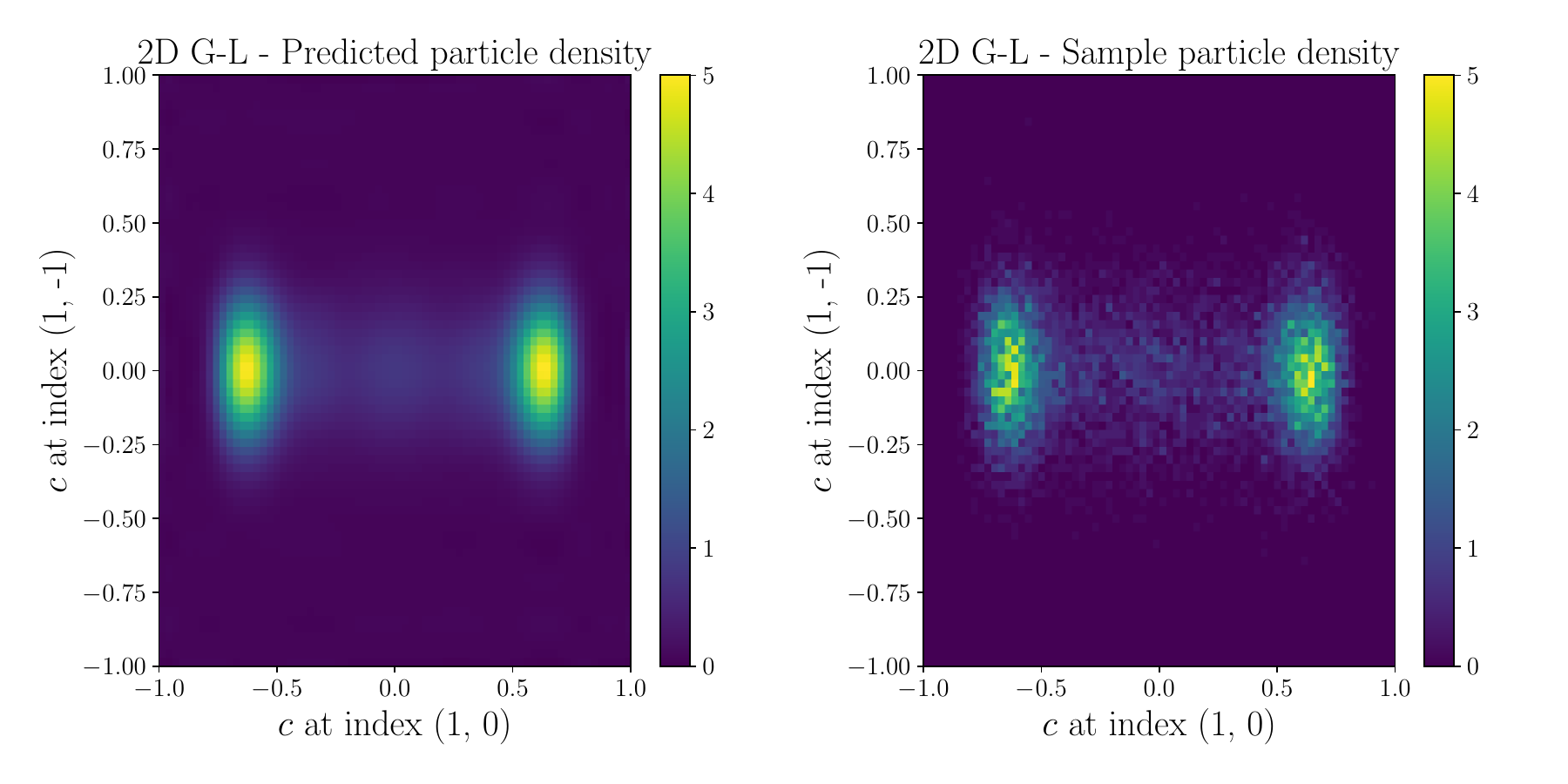}
  \caption{2D Ginzburg-Landau model. Plots of the marginal distribution of \(C \sim p(c)\) at \((c_{15, 5}, c_{8, 4})\), \((c_{2, 2}, c_{1, 1})\) and \((c_{1, 0}, c_{1, -1})\). For illustration purposes, a scaling is performed so that the marginal distribution lies in \([-1,1]^2\).}
  \label{Fig: 2D marginal}
\end{figure}

\subsection{2D lattice models}\label{sec: 2D lattice model result}

This subsection tests \Cref{alg:TTN density estimation internal} on strongly coupled 2D lattice models. We first consider the strongly coupled 2D Ginzburg-Landau model. We use the parameter setting as in \Cref{sec: 2D numerical rank}. The probability distribution function is 
\[
\begin{aligned}
    p(x_{(1, 1)}, \ldots, x_{(m, m)}) 
    =&\exp\Big(-\frac{\alpha_1}{2}\sum_{j}\sum_{i \sim i'}(x_{(i, j)}   - x_{(i', j)})^2  \\ & \,\,\,\,\,\,\,\,\,\,\,\,\,\,\,-\frac{\alpha_2}{2}\sum_{i}\sum_{j \sim j'}(x_{(i, j)} - x_{(i, j')})^2
- \frac{\lambda}{2}\sum_{i,j = 1}^{m}\left(1-x_{(i,j)}^2\right)^2\Big),
\end{aligned}
\] 
where \(m = 8\) with \(d = m^2 = 64\), \(\alpha_1 = 20\), \(\alpha_2 = 0.6\) and \(\lambda = 1\). From the discussion in \Cref{sec: 2D numerical rank}, we see that the numerical rank of this model makes it challenging for an FHT ansatz to perform density estimation. However, the numerical rank under \(p(c)\) is relatively small, which shows that the FHT-W ansatz is a suitable numerical candidate. 

\begin{figure}[h]
  \centering
  \includegraphics[width = \linewidth]{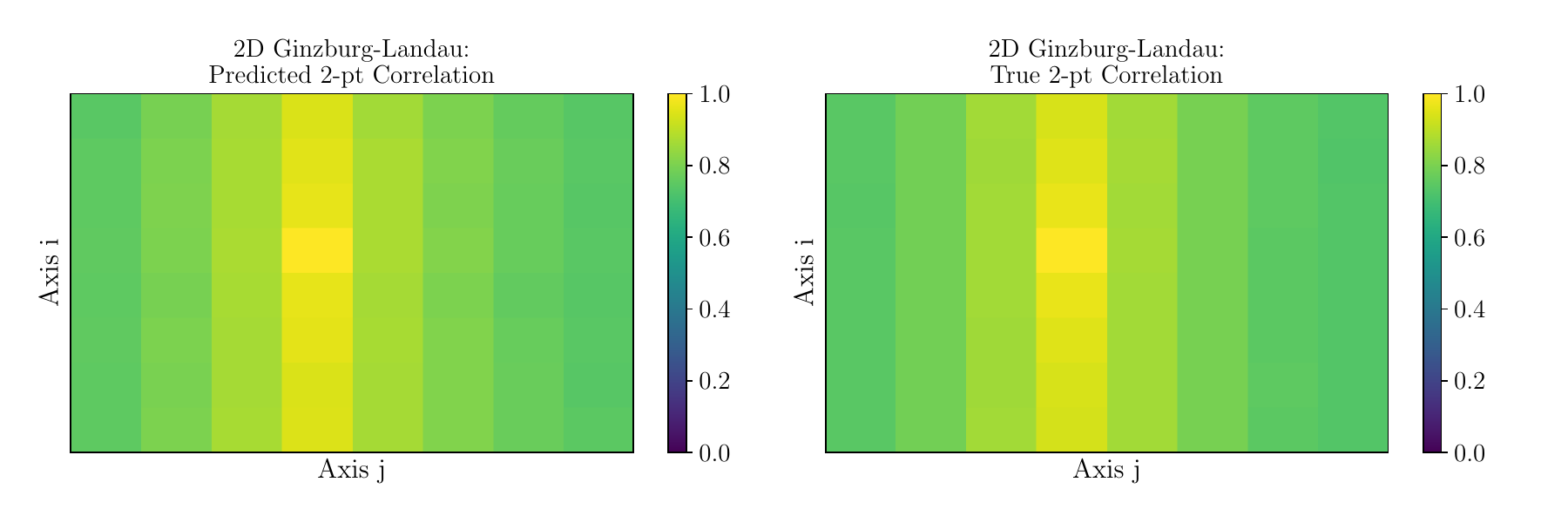}
  \caption{2D Ginzburg-Landau model. Plot of the two-point correlation function predicted by the FHT-W ansatz.}
  \label{Fig: 2D_GL_2_pt_correlation_64}
  \centering
  \includegraphics[width = \linewidth]{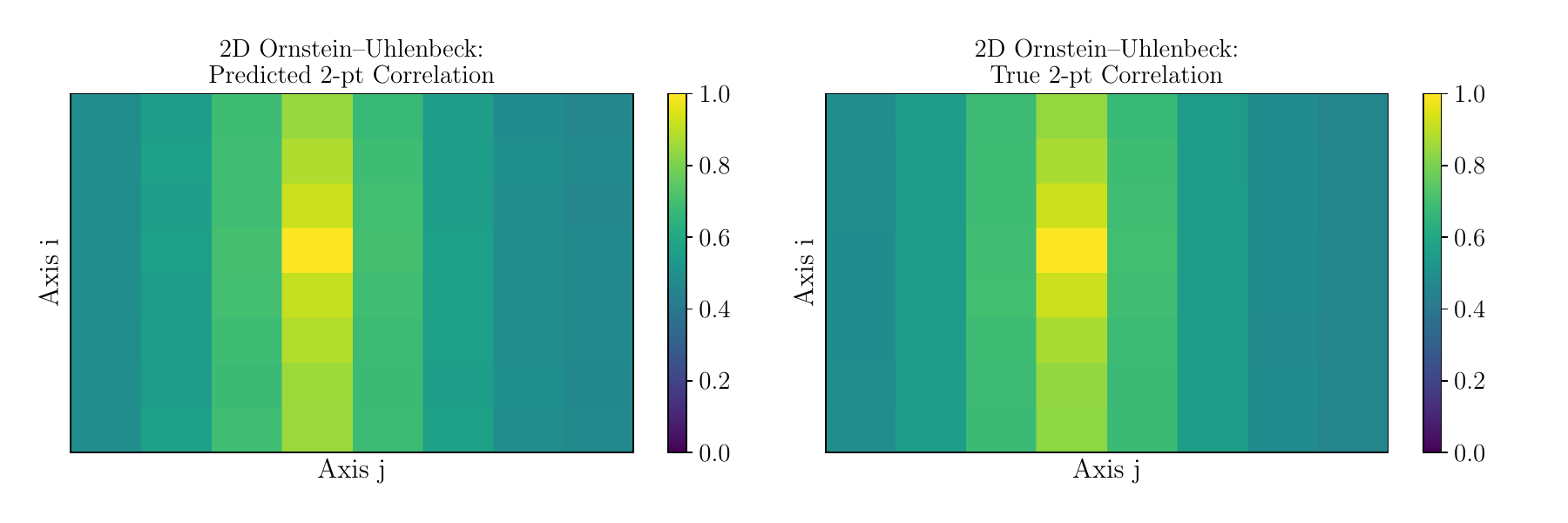}
  \caption{2D Ornstein–Uhlenbeck model. Plot of the two-point correlation function predicted by the FHT-W ansatz.}
  \label{Fig: 2D_OU_2_pt_correlation_64}
\end{figure}

The procedure to perform density estimation of 2D G-L models follows that of the 1D case. We first perform iterative wavelet transformations on \(N = 12000\) samples. We then perform density estimation on the wavelet transformed samples using the FHT-W ansatz architecture as specified in \Cref{sec: 2D lattice model main formulation}. Following the numerical experiment in \Cref{sec: 2D numerical rank}, we choose the maximal internal dimension to be \(r = 20\). 
On each variable, we use a Legendre polynomial basis with a maximal degree of \(q = 25\). By the construction in \Cref{sec: 2D lattice model main formulation}, we see that the 2D case also has a variable structure \(c = (c_{2L-1}, \ldots, c_{0}, c_{-1})\) with \(c_{l} = (c_{1, l}, \ldots, c_{2^l, l})\) for \(l = 1, \ldots, 2L-1\).

\Cref{Fig: 2D marginal} shows that the FHT-W ansatz accurately captures the marginal distribution of \(C \sim p(c)\) at \((c_{15, 5}, c_{8, 4})\), \((c_{2, 2}, c_{1, 1})\) and \((c_{1, 0}, c_{1, -1})\), which shows that the FHT-W ansatz is capable of capturing the interaction on fine scales and coarse scales.

To give a challenging test to the proposed method, we use the obtained FHT-W ansatz to predict the value of the two-point correlation function,
\[
f(i, j) =\mathrm{Corr}_{X \sim p(x)}\left(X_{(i,j)}, X_{(4, 4)}\right),
\]
which is done by using the obtained ansatz to perform repeated observable estimation tasks. \Cref{Fig: 2D_GL_2_pt_correlation_64} shows that the proposed FHT-W ansatz closely matches the ground truth. From the plot, we see that the distribution \(X = (X_{(i,j)}) \sim p(x)\) has a peculiar correlation structure of having a strong coupling on the \(i\) index but a weak coupling at the \(j\) index, which is one of the underlying reasons for the large numerical rank in \Cref{sec: 2D numerical rank}. 

We repeat the experiment on the 2D Ornstein–Uhlenbeck model, where the model parameter likewise follows from that of \Cref{sec: 2D numerical rank}.
We perform \Cref{alg:TTN density estimation internal} on \(N = 12000\) samples of \(p\). In \Cref{Fig: 2D_OU_2_pt_correlation_64}, we see that the FHT-W ansatz can likewise accurately model the two-point correlation function of the 2D O-U model.

\section{Conclusion}

This paper introduces a functional hierarchical tensor network under a wavelet basis for density estimation on lattice models. The approach combines a specific tree-based functional tensor network with wavelet transformation. The algorithm is applied successfully to the discretized Ornstein–Uhlenbeck model and Ginzburg-Landau model in 1D and 2D. With appropriate numerical treatment, this approach has the potential for modeling high-dimensional lattice models with complex coupling structures. An open question is whether one can use the FHT-W ansatz to efficiently solve the backward Kolmogorov equation and the Hamilton-Jacobi equation on lattice models.

\bibliographystyle{siamplain} 
\bibliography{references}

\end{document}